\documentclass[a4paper]{amsart}

\usepackage{enumerate}
\usepackage{todonotes}
\usepackage{appendix}
\usepackage{amssymb}

\newtheorem{thm}{Theorem}[section]
\newtheorem{lemma}[thm]{Lemma}

\newtheorem{remark}[thm]{Remark}
\newtheorem{proposition}[thm]{Proposition}
\newtheorem{example}[thm]{Example}

\newcommand{\N}{{\mathbb N}}

\newcommand{\R}{{\mathbb R}}

\newcommand{\cH}{{\mathcal H}}

\newcommand{\rn}{{\mathbb R^n}}

\allowdisplaybreaks

\DeclareMathOperator{\inte}{int}

\DeclareMathOperator*{\esssup}{ess\,sup}
\DeclareMathOperator*{\essinf}{ess\,inf}

\DeclareMathOperator{\dom}{dom}

\newcommand{\sprod}[2]{\langle #1,#2 \rangle} 

\newcommand{\Tu}{u^{\scriptscriptstyle K}}

\newcommand{\smallK}{{\scriptscriptstyle K}}
\newcommand{\smallL}{{\scriptscriptstyle L}}

\newcommand{\cutoff}[2]{{T}_{#1,#2}}

\numberwithin{equation}{section}

\begin{document}

\title{
Anisotropic symmetrization, convex bodies, and isoperimetric inequalities
}
\author[Gabriele Bianchi, Andrea Cianchi and Paolo Gronchi]
{Gabriele Bianchi, Andrea Cianchi and Paolo Gronchi}
\address{Dipartimento di Matematica e Informatica ``U.Dini", Universit\`a di Firenze, Viale Morgagni 67/A, Firenze, Italy I-50134} \email{gabriele.bianchi@unifi.it, andrea.cianchi@unifi.it}

\begin{abstract} This work is concerned with a P\'olya-Szeg\"o type inequality for anisotropic functionals of Sobolev functions. The relevant inequality entails a double-symmetrization involving both trial functions and functionals. A new approach  that uncovers geometric aspects of the inequality is proposed. It relies upon anisotropic isoperimetric inequalities, fine properties of Sobolev functions, and results from the Brunn-Minkowski theory of convex bodies. Importantly, unlike  previously available  proofs, the one offered in this paper does not require approximation arguments and hence allows for a characterization of extremal  functions.
\end{abstract}

\maketitle

\section{Introduction}

A classical functional inequality, usually called P\'olya-Szeg\"o principle in the literature, asserts that a convex Dirichlet type integral of a Sobolev function  $u: \R^n \to \R$, which decays in a suitable weak sense near infinity, does not increase under radially decreasing symmetrization. More specifically,
\begin{align}
    \label{PS}
    \int_{\R^n}A(|\nabla u^\bigstar|)\, dx \leq  \int_{\R^n}A(|\nabla u|)\, dx
\end{align}
 for any Young function $A: [0, \infty) \to [0, \infty]$.
Here, $\nabla u$ denotes the weak gradient of $u$, $|\nabla u|$ its Euclidean norm, and $u^\bigstar$ the radially decreasing symmetral of $u$, i.e. the function equimeasurable with $u$ whose level sets are  balls centered at the origin. The definition of Young function, as well as other notions appearing in this introduction, are recalled in the next section. We shall here content ourselves with
 mentioning that the choice $A(t)=t^p$ is admissible in \eqref{PS} for any $p \geq 1$.
The inequality \eqref{PS} is of critical use in the proof of a number of results concerning sharp constants and optimal shapes of geometric objects in mathematical analysis, differential geometry, and mathematical physics.

Although diverse approaches to \eqref{PS} are available, a customary proof relies upon a combination of the coarea formula with the classical isoperimetric inequality -- see e.g. \cite{BrothersZiemer, CiFuAppAn, Talenti}. The appeal of this method stems from shedding light on the geometric flavor of the inequality \eqref{PS}. In fact, when $A(t)=t$, this inequality can be regarded as a functional form of the isoperimetric inequality in $\R^n$.

Variants and generalizations of the inequality \eqref{PS}, involving different integral functionals, broader classes of trial functions, alternate  symmetrizations are well known and are spread in a vast literature. A unified approach embracing a broad class of symmetrizations can be found in \cite{BiGaGrKi2022}. Reference monographs on this topic are \cite{Ka, Ke, Bae}.
Unconventional strengthened versions of the P\'olya-Szeg\"o principle,  concerning energy functionals invariant under affine transformations and hence called affine P\'olya-Szeg\"o principles, were introduced in  \cite{Lyz, Zhang}.

A popular extension of \eqref{PS} allows for  norms $H(\nabla u)$ of $\nabla u$ more general than just the Euclidean one. It tells us that
\begin{align}
    \label{PSH}
    \int_{\R^n}A(H(\nabla u^\circ))\, dx \leq  \int_{\R^n}A(H(\nabla u))\, dx,
\end{align}
where $u^\circ$ denotes the symmetral of $u$ with respect to the dual norm $H_0$ to $H$.  Namely, $u^\circ$ is the function equimeasurable with $u$, whose super-level sets are concentric balls in the metric induced by $H_0$,  i.e., dilates of the unit ball $\{ H_0\leq 1\}$. These balls are usually called the Wulff shapes associated with $H$, after the name of T.~Wulff who introduced anisotropic perimeters in \cite{Wu}.
The inequality \eqref{PSH}
was established in \cite{AFLT}, under the name of convex symmetrization inequality, via arguments  along the same lines as those for \eqref{PS}, save that the standard isoperimetric inequality is replaced with an anisotropic isoperimetric inequality dictated by the norm $H$.

The present paper deals with an even more general, fully anisotropic inequality, where the integrand in the   functional is an $n$-dimensional Young functions  $\Phi : \R^n \to [0, \infty]$. Hence, the functional depends on the whole gradient $\nabla u$ and not just on its norm. The  Dirichlet type integral considered have the form:
\begin{align}
    \label{integral}
    \int_{\R^n} \Phi (\nabla u)\, dx.
\end{align}
Loosely speaking, the inequality in question ensures that the functional \eqref{integral} does not increase if both $u$ and $\Phi$ are simultaneously properly symmetrized with respect to any given convex body $K$. The symmetral $u^\smallK$ of $u$ is just the function equimeasurable with $u$, whose super-level sets are homothetic to $K$. The symmetrization of the integrand $\Phi$ is accomplished through a less straightforward three-step process. This process calls into play, in a sequence, the operation of Young conjugation (also called Legendre transform), a symmetrization with respect to $K$, and a Young conjugation again. Altogether, this yields a new $n$-dimensional Young function that will be denoted by $\Phi_{\bullet \smallK \bullet}$. The anisotropic symmetrization inequality then reads:
\begin{align}
    \label{PSaniso}
 \int_{\R^n} \Phi_{\bullet \smallK \bullet} (\nabla u^\smallK)\, dx \leq     \int_{\R^n} \Phi (\nabla u)\, dx.
\end{align}

This is a kind of universal P\'olya-Szeg\"o principle, which encompasses various special instances, including \eqref{PS} and \eqref{PSH}.
For $K$ equal to a ball, the inequality \eqref{PSaniso} was stated, without a proof, by V.~Klimov in his paper \cite{Kli74}. Over the years, Klimov has authored several papers around this inequality and its applications, but,
as far as we know, he   never published a full proof, although some hints are given in \cite{Kli74}. In \cite{Kli99} he established a variant of \eqref{PSaniso}, where symmetrization with respect to a convex set is replaced with Steiner symmetrization about a hyperplane. Thanks to the latter result, a proof of \eqref{PSaniso}, when $K$ is a ball, was accomplished in \cite{CiCPDE} via an approximation process by sequences of Steiner symmetrizations. In the paper \cite{VS2006}, the inequality \eqref{PSaniso} is proved for a general convex body $K$. The argument of \cite{VS2006} again rests upon  approximations, which exploit  sequences of polarizations.

  The inequality \eqref{PSaniso} is the point of departure in the approach to Sobolev type inequalities in anisotropic Orlicz-Sobolev spaces. Early results in this direction are contained in \cite{Kli74}; sharp embeddings are the subject of \cite{cianchi_pacific, cianchi_ibero}. These embeddings have a role in the existence and regularity theory of solutions to anisotropic elliptic equations and variational problems -- see e.g. \cite{alberico,  ADF, ACCZ, barletta, CiAIHP, CiCPDE}.

The purpose of this work is to offer an alternate, direct proof of the inequality \eqref{PSaniso}, which avoids approximation arguments via partial symmetrizations. The proof to be presented provides geometric insight on \eqref{PSaniso}, which, by contrast, does not emerge form the previously available approaches. It rests on the anisotropic isoperimetric inequality mentioned above  and techniques which pertain to the Brunn-Minkowski theory of convex bodies. Results from geometric measure theory and fine properties of Sobolev functions are also exploited, since our proof does not require regularizations of  $u$  and
applies directly to Sobolev functions.

The absence of any approximation argument is a major advantage of the  new proof, which  enables us to characterize those functions $u$ for which equality holds in \eqref{PSaniso}.
The question of the cases of equality in symmetrization inequalities is a delicate issue. It is well known that, if equality holds in \eqref{PS}, then the super-level sets of $u$ are balls and $|\nabla u|$ is constant on their boundaries. A parallel result holds with regard to \eqref{PSH}, save that Euclidean balls have to be replaced with balls in the metric of $H_0$, and $|\nabla u|$ with $H(\nabla u)$. Through a close inspection of the steps in our proof of \eqref{PSaniso} we can  provide information in a similar spirit about functions which render the inequality \eqref{PSaniso}
an equality. Due to the freedom in the choice of  $K$ and $\Phi$, and the presence of a double symmetrization, the characterization of the extremals in this inequality is inevitably subtler.

Recall that full symmetry of extremals in \eqref{PS} and \eqref{PSH}, namely the fact that their level sets be concentric, is not guaranteed without additional assumptions on $A$ and $u$. The correct minimal assumptions to be imposed were exhibited in \cite{BrothersZiemer} for the inequality \eqref{PS}, and extended to the inequality \eqref{PSH} in \cite{EsTr} and \cite{FeVo}. They amount to requiring that $A$ be strictly convex and that the set of critical points of $u^\bigstar$ (or $u^\circ$) have measure zero. Of course, a symmetry result for extremals in \eqref{PSaniso} in the same direction is an interesting problem, which, however, falls beyond the scope of the present paper. In fact, as will be apparent from the examples that will be produced, it is not  clear how a possible statement of a result in this direction should read.

\bigskip

\emph{In memory of Paolo Gronchi.} Paolo passed away on July 4th, 2024, shortly before the final version of this paper was completed. He was an excellent mathematician and has been a close companion of ours for many years.  It has been a pleasure to work and spend time with him. We have always enjoyed his constant good humor and positivity.

\section{Notation and background}

In this section we introduce some notation and recall basic properties from the theory of convexity and of Sobolev type spaces.

\subsection{Gauge functions and convex bodies}
A function $H : \R^n \to [0, \infty)$ will be called a \emph{gauge function} if it is convex, positively homogeneous of degree $1$,  and vanishes only at $0$. Consequently, for any gauge function $H$ there exists positive constants $c_1$ and $c_2$ such that
\begin{align}
    \label{equivH}
c_1 |\xi| \leq H(\xi) \leq c_2 |\xi| \qquad \text{for $\xi \in \rn$.}
\end{align}
 The dual of a gauge function $H$, denoted by $H_0$, is also a gauge function and is defined as
\begin{align}
    \label{H0}
    H_0(\xi)= \max_{\R^n \ni \eta \neq 0}\frac {\xi \cdot \eta}{H(\eta)} \qquad \text{for $\xi \in \R^n$.}
\end{align}
Notice that
$$(H_0)_0=H.$$
A \emph{convex body} is  a closed bounded convex set in $\R^n$. Let $L$ be a convex body such that $0\in \inte (L)$,
where
$\inte (L)$ denotes
the set of the interior points of  $L$.
\\ The \emph{polar body} of $L$ is the convex body $L^\circ$ defined as
$$L^\circ=\{x\in\R^n : \sprod{x}{y}\leq1\text{ for all $y\in L$}\}.$$
The \emph{gauge function of $L$} is the function $H_0^\smallL : \R^n \to [0, \infty)$ defined as
$$ H_0^\smallL(x)=\min\{\lambda\geq0\ : x\in\lambda L\} \qquad \text{ for $x\in \R^n$.}$$
The function $H_0^\smallL$ is a gauge function in the sense specified above.
\\ The  \emph{support function} $h_L : \R^n \to [0, \infty)$ of $L$ is defined by
$$h_L(\xi)=\max\{\sprod{\xi}{\eta}: \eta\in L\} \qquad \text{for $\xi\in\R^n$.}$$
One has that $h_L=H^\smallL$, where $H^\smallL$ denotes the dual of $H_0^\smallL$.
Moreover, $H_0^\smallL$ agrees with the support function of $L^\circ$, whereas $h_L$ agrees with the gauge function of $L^\circ$ \cite[Cor.~14.5.1]{Rockafellar}.
Also,
\begin{align}
L&=\{x : H_0^\smallL(x)\leq 1\}\label{L=h0Llq1}\\
\intertext{and}
L^\circ&=\{\xi : h_L(\xi)\leq1\},\label{polarL=hLlq1}
\end{align}
see \cite[Cor.~15.1.2]{Rockafellar}.\\
If $\xi \in \rn \setminus \{0\}$ is such that $\nabla H_0(\xi)$ exists, then
\begin{equation}
\label{nablaH}
    H(\nabla H_0(\xi))=1,
\end{equation}
  see e.g. \cite[Equation (3.12)]{CianchiSalani}.
  \\
The anisotropic isoperimetric inequality associated with $h_L$ tells us that
\begin{equation}\label{may9}
\int_{\partial^*E} h_L(\nu^E)\, d \cH^{n-1} \geq n |E|^{\frac{n-1}{n}} |L|^{\frac 1n}
\end{equation}
for every set  $E\subset \R^n$ of finite perimeter. Here, $\partial ^*E$ denotes the reduced boundary of $E$, and  $\nu^E$ stands for the geometric measure theoretical outer unit normal to $E$. Moreover, equality holds in \eqref{may9} if and only if $E$ is a dilate of $L$,   up to sets of Lebesgue measure zero and up to translations -- see \cite{FonMul1991}.

\subsection{Young functions}
 A function $$A : [0, \infty) \to [0, \infty]$$ is said to be a \emph{Young function} if it is convex, left-continuous, $A(0)=0$ and non-constant in $(0,\infty)$. The \emph{Young conjugate} of $A$, which will be denoted by $A_\bullet$, is the Young function obeying:
 \begin{align}
     \label{Aconj}
 A_\bullet (t) = \sup\{st - A(s): s\geq 0\} \qquad \text{for $t\geq 0$.}
 \end{align}
 One has that
 $$A=(A_\bullet)_\bullet$$
 for every Young function $A$. Therefore, Young conjugation is an involution.
 \\
Besides plain power functions $t^p$, with $p\geq 1$, basic instances of Young functions are
\begin{align}
    \label{powerlog}
    A(t) = t^p (\log (c+t))^q \qquad \text{for $t \geq 0$,}
\end{align}
where the positive constant $c$ is large enough for $A$ to be convex, and
\begin{align}
    \label{exp}
    A(t) = e^{t^\alpha}-1 \qquad \text{for $t \geq 0$,}
\end{align}
with $\alpha >0$.

\subsection{ $n$-dimensional Young functions}

We call an \emph{$n$-dimensional Young function} a function
 $$\Phi:\R^n\to [0, \infty]$$
 which is
 convex,  lower semicontinuous,  finite in a neighborhood of $0$, and such that
\begin{align}\label{hpPhi}
   \text{ $\Phi(0)=0$ \,\, and \,\, $\lim_{|\xi|\to\infty}\Phi(\xi)=+\infty$.}
\end{align}
The domain of $\Phi$ is denoted by $\dom \Phi$ and defined as
$$\dom\Phi=\{\xi\in\R^n : \Phi(\xi)<\infty\}.$$
Notice that for any $n$-dimensional Young function  there exist positive constants $c$ and $r$ such that
\begin{align}
    \label{2024-30prima}
    \Phi(\xi) \geq c|\xi| \qquad \text{if $|\xi|\geq r$.}
\end{align}
 The subgradient of $\Phi$ at a point $\xi \in \R^n$ will be denoted by $\partial \Phi (\xi)$. Recall that $\partial \Phi (\xi)\subset \R^n$ and, if $\Phi$ is differentiable at $\xi$, then $\partial \Phi (\xi)= \{\nabla \Phi (\xi)\}$.
\\
For the \emph{Young conjugate} of an $n$-dimensional Young function $\Phi$ we employ, with a slight abuse,  the same  notation as in \eqref{Aconj}, namely $\Phi_\bullet$. It is defined as
$$\Phi_\bullet (\xi) = \sup \{\sprod{\xi}{\eta} - \Phi(\eta): \eta \in \R^n\} \qquad \text{for $\xi \in \R^n$,}$$
where $\sprod{\cdot}{\cdot}$ denotes scalar product in $\R^n$.   The function  $\Phi_\bullet$ inherits the properties of $\Phi$, and hence it is also an  $n$-dimensional Young function \cite[Corollary 14.2.2]{Rockafellar}.
\\
The operation of Young conjugation in the class of $n$-dimensional Young functions is an involution, since
\begin{equation}
    \label{2024-38}
\Phi=\big(\Phi_\bullet\big)_\bullet.
\end{equation}
Observe that, if $A$ is a Young function, then the function $\Phi$ given by
$$\Phi (\xi) = A(|\xi|) \qquad  \text{for $\xi \in \R^n$}$$
 is an $n$-dimensional Young function.
 \\ Genuinely anisotropic instances of $n$-dimensional Young functions have the form
\begin{equation}\label{Phi=Ai}
\Phi (\xi) = \sum _{i=1}^n A_i(|\xi_i|) \qquad \hbox{for $\xi \in \R^n$,}
\end{equation}
where $A_i$ are Young functions and $\xi = (\xi_1, \dots \xi_n)$.
A customary example of functions of this kind is
\begin{equation}\label{Phi=pi}
\Phi (\xi) = \sum _{i=1}^n |\xi_i|^{p_i} \qquad \hbox{for $\xi \in \R^n$,}
\end{equation}
where  $1\leq p_i<\infty$, for  $i=1,\dots,n$.
\\
The example \eqref{Phi=Ai} can be further generalized
as
\begin{equation}\label{Phi=full}
\Phi (\xi) =  \sum _{k=1}^n A_k\Big(\Big|\sum_{i=1}^n \alpha _{ki}\xi_i\Big|\Big) \qquad \hbox{for $\xi \in \rn$,}
\end{equation}
where $A_k$ are Young functions, $m \in \mathbb N$, and   the matrix $(\alpha _{ki}) \in \mathbb R^{n\times n}$ is such that ${\rm det} (\alpha _{ik})\neq 0$. A possible instance, for $n=2$, is
\begin{equation}
\label{trud}
\Phi (\xi) = |\xi_1 -\xi_2|^p + |\xi_1|^q\log (c+ |\xi _1|)^\alpha \quad \hbox{for $\xi \in \mathbb R^2$,  }
\end{equation}
where either $q\geq 1$ and $\alpha >0$, or $q=1$ and $\alpha >0$, the exponent  $p\geq 1$, and $c$ is a sufficiently large constant for $\Phi$ to be convex. Another example   amounts to the function
\begin{equation}
\label{trud1}
\Phi (\xi) = |\xi_1 + 3 \xi_2|^p + e^{|2\xi_1-\xi_2|^\beta} -1  \quad \hbox{for $\xi \in \mathbb R^2$,  }
\end{equation}
with $p\geq 1$ and $\beta >1$.
\\ However, let us stress that there exist $n$-dimensional Young functions  which do not split as in \eqref{Phi=full} -- see e.g. \cite{ChNa}

\subsection{Sobolev functions}
We denote by $\mathcal M(\R^n)$ the set of real-valued measurable functions in $\R^n$. The \emph{distribution function}
$\mu : (\essinf u, \infty) \to [0, \infty]$ of a function $u \in {\mathcal{M}}(\R^n)$  is given by
\begin{align}\label{mu}
   \mu (t) = |\{x:u(x)>t\}|\qquad \text{for $t> \essinf u$.}
\end{align}
Define
$${\mathcal{M}_d}(\R^n) =\{u\in \mathcal M(\R^n): \mu (t)<\infty \,\,\text{for} \,\,t>\essinf u\}.$$
The set ${\mathcal{M}_d}(\R^n)$ can be regarded as the subset of those functions from $\mathcal M(\R^n)$ which decay to $\essinf u$ near infinity in the weakest possible sense in view of our applications.
\\ The decreasing
rearrangement $u^* : [0, \infty) \to [0, \infty]$ of a function $u\in {\mathcal{M}_d}(\R^n)$ is the generalized right-continuous inverse of $\mu$, which obeys
\begin{align}
    \label{u*}
    u^*(s)= \inf\left\{t\in \R:  \mu (t)\leq s\right\} \quad\text{for $s\geq 0$.}
\end{align}
Namely, $u^*$ is the non-increasing right-continuous function on $[0, \infty)$ such that $|\{s\in [0, \infty): u^*(s)>t\}|=|\{x\in \rn: u(x)>t\}|$ for $t>\essinf u$.
\\ The increasing rearrangement $u_*$ is defined analogously for any function $u$ such that $(-u) \in \mathcal M_d(\rn)$.
\\
A function $u\in \mathcal M(\rn)$ will be called \emph{quasi-convex} if
\begin{align}
    \label{quasic}
    \text{the set $\{u\geq t\}$ is convex for   $t \in (\essinf u, \esssup u)$.}
\end{align}
Given an $n$-dimensional Young function $\Phi$, we define
 the \emph{homogeneous Sobolev class}
$$V^{1,\Phi}(\R^n) = \bigg\{u\in W^{1,1}_{\rm loc}(\rn):   \int_{\R^n}\Phi(\nabla u)\, dx<\infty \bigg\},$$
and
set
$$V^{1,\Phi}_{\rm d}(\R^n)=
V^{1,\Phi}(\R^n)\cap\mathcal{M}_d(\R^n).$$
Notice that the sets $V^{1,\Phi}(\R^n)$ and $V^{1,\Phi}_{\rm d}(\R^n)$ are convex, but, because of the lack of homogeneity of $\Phi$, they need not be  linear spaces. In this connection, let us point out that, by contrast, in earlier contribution the notation $V^{1,\Phi}(\R^n)$ and $V^{1,\Phi}_{\rm d}(\R^n)$ was employed to denote the homogeneous Sobolev spaces associated with the function $\Phi$.
\\ When $\Phi (\xi)=|\xi|$, we shall simply write $V^{1,1}(\rn)$ and $V^{1,1}_{\rm d}(\rn)$ to denote the corresponding  Sobolev classes, which, in this case, agree with the respective Sobolev spaces.
\\ Given an $n$-dimensional Young function $\Phi$, a function $u\in V^{1,\Phi}_{\rm d}(\R^n)$, and $t> \essinf u$,
the function $\max \{u, t\}$ attains the constant value $t$ outside a set of finite measure. As a consequence of \eqref{2024-30prima} and of standard properties of truncations of Sobolev functions, we hence have that
\begin{align}
    \label{2024-31}
    \max \{u, t\}\in V^{1,1}_{\rm d}(\rn) \quad \text{if $t> \essinf u$}.
\end{align}
 A property of Sobolev functions ensures that, if  $u\in V^{1,\Phi}_{\rm d}(\R^n)$, then $u$ admits a representative such that
the set $\{u>t\}$ is of finite perimeter   and
\begin{align}\label{2024-6prima}
\partial^*\{u>t\} = \{u=t\} \qquad \text{up to sets of $\cH^{n-1}$ measure zero,}
\end{align}
for a.e. $t> \essinf u$ \cite{BrothersZiemer}. In what follows, by $u$ we always denote such a representative.
\\ The coarea formula ensures that, if  $u\in V^{1,1}_{\rm d}(\R^n)$, then
\begin{align}
    \label{coarea}
    \int_{\rn} f |\nabla u|\, dx = \int_{-\infty}^\infty\int_{\{u=t\}}f\, d\mathcal H^{n-1}\, dt
\end{align}
for every Borel function $f: \rn \to [0, \infty)$.

\subsection{Symmetrizations}

Let $K\subset\R^n$ be a convex body such that $0\in \inte (K)$.
 Given a measurable set $E\subset\R^n$, we  define $$E^\smallK=c K,$$ where $c\geq0$ is such that
$$|E|=|E^\smallK|.$$
 Here, $|\,\cdot\,|$ stands for  Lebesgue measure.
\\
 Given a function $u\in {\mathcal{M}_d}(\R^n)$ we define $\Tu$ as the symmetral of $u$ with respect to $K$. Namely, $\Tu : \R^n \to \R$  is the function obeying:
$$\{x\in\R^n : \Tu(x)\geq  t\}= \{x\in\R^n : u(x)\geq t\}^\smallK  \qquad \text{for   $t\in (\essinf u, \esssup u)$}.$$
%
Moreover,  if $\Phi$ is an $n$-dimensional Young function, we denote by $\Phi_\smallK : \R^n \to \R$ the function given by
$$-\Phi_\smallK(-\xi)=(-\Phi)^\smallK(\xi) \qquad \text{for $\xi\in \R^n$.}$$
Equivalently, $\Phi_\smallK$ is the symmetral of $\Phi$ such that
$$\{\xi\in\R^n : \Phi_\smallK(\xi)\leq t\}=- \{\xi\in\R^n : \Phi (\xi)\leq t\}^\smallK \qquad \text{for   $t\geq 0$.}$$
The latter symmetrization will be applied
in combination with the operation of Young conjugation. For simplicity of notation, given an  $n$-dimensional Young function $\Phi$, we set
\[\Phi_{\bullet \smallK}=\big(\Phi_\bullet\big)_\smallK\quad\quad\text{and}\quad\quad
\Phi_{\bullet \smallK \bullet}=\big(\big(\Phi_\bullet\big)_\smallK\big)_\bullet.
\]
Observe that, if $K$ is a ball, then $u^\smallK = u^\bigstar$, the radially decreasing symmetral of $u$, and $\Phi_\smallK = \Phi_\bigstar$, the radially increasing symmetral of $\Phi$.
\\
The operations of Young conjugation and symmetrization with respect to a convex set do not commute, even if $K$ is a ball.
Hence, $\Phi_{\bullet \smallK \bullet} \neq \Phi_\smallK$ in general. However, the  functions $\Phi_{\bullet \smallK \bullet}$ and $ \Phi_\smallK$
are equivalent, up to  constants multiplying their arguments. Namely, for every convex set $K$ as above, there exist constants $c_1$ and $c_2$, depending on $K$, such that
\begin{align}
\label{equivalence}
\Phi_\smallK (c_1\xi) \leq \Phi_{\bullet \smallK \bullet} \leq \Phi_\smallK (c_2\xi) \qquad \text{for $\xi \in \rn$,
}
\end{align}
for every $n$-dimensional Young function $\Phi$. Equation \eqref{equivalence}
is established in
\cite[Lemma 7]{Kli76} when $K$ is a ball.
On the other hand, analogously to the formula \eqref{uK_and_gauge} below, one has that
\begin{align*}
    \Phi_\smallK (\xi) = \Phi_*(\kappa H_0(\xi)^n) \qquad \text{for $\xi \in \rn$,
}
\end{align*}
where $\Phi_*$ is the increasing rearrangement of $\Phi$, $H_0$ is the gauge function of $K$, and $\kappa = |K|$. Hence, thanks to the property
 \eqref{equivH} applied to $H_0$, the functions $\Phi_\smallK$ and $\Phi_\smallL$ are equivalent for any couple of convex bodies $K$ and $L$. Equation \eqref{equivalence} for balls thus implies it for an arbitrary $K$.

\section{Main results}

The anisotropic symmetrization principle which is the subject of this paper
is the content of Theorem \ref{t:klimov} below. In its statement, and in what follows,
$K$
denotes a convex body such that $0\in \inte (K)$.

\begin{thm}\label{t:klimov}
Let $\Phi$ be an $n$-dimensional Young function.
Assume that $u\in V^{1,\Phi}_{\rm d}(\R^n)$.  Then $\Tu\in V^{1,\Phi_{\bullet \smallK \bullet}}_{\rm d}(\R^n)$ and
\begin{equation}\label{ineq_klimov}
\int_{\R^n}\Phi_{\bullet \smallK \bullet}(\nabla \Tu)\, dx\leq \int_{\R^n} {\Phi}(\nabla u)\, dx.
\end{equation}
\end{thm}

As mentioned above, one major benefit from the proof of the inequality \eqref{ineq_klimov} that will be offered is in the possibility of characterizing those functions $u$ which render this inequality  an equality. In the next result necessary conditions on $u$ are exhibited for equality to hold in  \eqref{ineq_klimov}. The subsequent Theorem \ref{t:klimov-suff} asserts that the relevant conditions are also sufficient.
\begin{thm}\label{t:klimov-nec-new}
Let $\Phi$ be an $n$-dimensional Young function such that
\begin{align}
    \label{2024-6}
    0<\Phi (x) < \infty \qquad \text{for $x \neq 0$,}
\end{align}
let $u\in  V^{1,\Phi}_{\rm d}(\R^n)$  and assume that equality holds in the inequality \eqref{ineq_klimov}.  Then  $u$ equals  a quasi-convex function a.e. in $\rn$.  Moreover:
\begin{enumerate}[(i)]
\item
For a.e. $t\in(\essinf u, \esssup u)$ there exist   $s_t\geq0$, $a_t>0$, and $x_t\in\R^n$ such that:
\begin{enumerate}[(a)]
\item  \begin{align}
    \label{2024-30}
    \text{$\inte (\{\Phi_\bullet\leq s_t\})\neq \emptyset$; }
\end{align}
\item \begin{equation}\label{t:klimov-suff-a}
\{u\geq t\}=-a_t\{\Phi_\bullet\leq s_t\}+x_t \qquad \text{  up to a set of Lebesgue measure zero;}
 \end{equation}
\item For $\cH^{n-1}$-a.e. $x\in\{u=t\}$ there exists  $\xi\in\{\Phi_\bullet = s_t\}$ such that
\begin{equation} \label{t:klimov-suff-c}
\nabla u(x)\in \partial\, \Phi_\bullet (\xi);
\end{equation}
\item For $\cH^{n-1}$-a.e. $x\in\{\Tu=t\}$ there exists  $\xi\in\{\Phi_{\bullet\smallK} = s_t\}$ such that
\begin{equation} \label{t:klimov-suff-d}
\nabla \Tu(x)\in \partial\, \Phi_{\bullet\smallK} (\xi);
\end{equation}
\item  If $s_t>0$, then
\begin{equation}\label{t:klimov-suff-b}
\int_{\{u=t\}} \frac{d \cH^{n-1}}{|\nabla u|}\,= \int_{\{\Tu=t\}} \frac{d \cH^{n-1}}{|\nabla \Tu|}.
\end{equation}
\end{enumerate}
\item Assume, in addition, that
 \begin{align}
     \label{strict}
     \text{$\Phi$ is strictly convex.}
 \end{align}
Then  $\Phi_\bullet$ is differentiable  in $\inte(\dom\Phi_\bullet)$
and \eqref{t:klimov-suff-c} can be replaced with
\begin{equation} \label{equality_necessity_grad1}
\nabla u(x)=\nabla \Phi_\bullet (\xi)
\end{equation}
for some $\xi \in \{\Phi_\bullet =s_t\} \cap \inte(\dom\Phi_\bullet)$.
\\
  In particular, if
\begin{align}
    \label{diff}
    \text{$\Phi$ is differentiable,}
\end{align}
then also $\Phi_\bullet$ is strictly convex  in $\inte(\dom\Phi_\bullet)$ and there exists a unique $\xi \in \{\Phi_\bullet =s_t\} \cap \inte(\dom\Phi_\bullet)$ fulfilling \eqref{equality_necessity_grad1}.
\end{enumerate}
\end{thm}

In Examples \ref{classical}
and \ref{convexsym} below, the classical P\'olya-Szeg\"o principle \eqref{PS} and  its extension to non-Euclidean norms are recovered from \eqref{ineq_klimov}. The well known information about their extremals is also deduced via Theorem \ref{t:klimov-nec-new}.

\begin{example}\label{classical}{\rm
Assume that
$$\text{$K$ is an Euclidean ball, centered at $0$}$$
and
$$\Phi (\xi) = A(|\xi|) \qquad \text{for $\xi \in \R^n$,}$$
for some classical Young function $A$. Since $\Phi$ is radially symmetric,  $\Phi_\bullet$ is also radially symmetric. Therefore,
$\Phi_{\bullet\smallK} = \Phi_\bullet$, inasmuch as $K$ is a ball.
Altogether,
$$\Phi_{\bullet\smallK \bullet}(\xi) = A(|\xi|) \qquad \text{for $\xi \in \R^n$.}$$
Moreover, if $u\in V^{1,\Phi}_d(\R^n)$, then
$$u^K = u^\bigstar ,$$
the  radially decreasing rearrangement of $u$. Therefore, the inequality \eqref{ineq_klimov} reproduces the classical P\'olya-Szeg\"o principle \eqref{PS}.
\\ Assume, in addition, that
\begin{align}
    \label{2024-83}
    0<A(t)<\infty \qquad \text{for $t>0$,}
\end{align}
and that equality holds in \eqref{PS} for some function $u\in V^{1,\Phi}_d(\R^n)$.
Then equation \eqref{t:klimov-suff-a} of Theorem \ref{t:klimov-nec-new} enable  us to deduce that
$$\text{$\{u\geq t\}$ is a ball}$$
for $t\in (\essinf u, \esssup u)$.
Moreover, if
\begin{align}
    \label{2024-84}
   \text{$A$ is strictly convex,}
\end{align}
then equation \eqref{equality_necessity_grad1} implies that
$$\text{$|\nabla u|$ is constant $\mathcal H^{n-1}$ -- a.e. on $\{u=t\}$}$$
for a.e. $t \in (\essinf u, \esssup u)$.
}
\end{example}

\begin{example}\label{convexsym}{\rm
Let $K$ be any convex body such that $0\in \inte (K)$ and let $H_0$ be its gauge function.
Assume that $A$ is a Young function and
$$\Phi (\xi) = A(H(\xi)) \qquad \text{for $\xi \in \R^n$,}$$
where $H: \R^n \to [0, \infty)$ is the dual of $H_0$.
One has that
\begin{align}
    \label{2024-80}
    \Phi_{\bullet\smallK \bullet} (\xi) = A(H(\xi)) \qquad \text{for $\xi \in \R^n$.}
\end{align}
  Indeed, \begin{align}\label{2024-200}
\Phi_\bullet(\xi)=A_\bullet(H_0(\xi))
  \end{align} -- see Lemma \ref{convesymm}, Section \ref{techyoung}. Thus,
 all the sub-level sets of $\Phi_\bullet$ are homothetic to $K$,  whence  $\Phi_{\bullet\smallK}=\Phi_\bullet$. The involution property of Young conjugation then yields $\Phi_{\bullet\smallK \bullet}=\Phi$, namely \eqref{2024-80}.
\\ On the other hand, if $u\in V^{1,1}_d(\R^n)$, then
$$u^K = u^\circ,$$
where
\begin{equation*}
u^\circ(x) = u^*(\kappa H_0(x)^n) \qquad \text{for $x\in \R^n$,}
\end{equation*}
and $\kappa = |K|$.
\\ As a consequence,
 the inequality \eqref{ineq_klimov} recovers  \eqref{PSH}.
\\ Assume now that \eqref{2024-83} is in force. Then, from  equation \eqref{t:klimov-suff-a} of Theorem \ref{t:klimov-nec-new} we can infer that
$$\text{
$\{u\geq t\}$ is a ball in the metric induced by $H_0$}$$
for $t\in (\essinf u, \esssup u)$.
 \\ If $A$ fulfils the stronger assumption \eqref{2024-84}, then equations \eqref{equality_necessity_grad1} and \eqref{2024-200} entails that
 $$\text{$\nabla u (x)=  A'_\bullet(H_0(x))\nabla H_0 (x)$ \,\, for $\mathcal H^{n-1}$ -- a.e.  $x\in\{u=t\}$}$$
 for a.e $t \in  (\essinf u, \esssup u)$.
\\
 Since $H$ is positively homogeneous of degree $1$ and satisfies property \eqref{nablaH}, we hence deduce that
 $$\text{
 $H(\nabla u)$ is constant $\mathcal H^{n-1}$ -- a.e. on $\{u=t\}$}$$
 for a.e. $t \in (\essinf u, \esssup u)$.}
\end{example}

\begin{thm}\label{t:klimov-suff}
Assume that  $\Phi$ is as in Theorem \ref{t:klimov} and fulfills the condition \eqref{2024-6}.
Let $u\in V^{1,\Phi}_{\rm d}(\R^n)$.
Assume that, for a.e.  $t\in(\essinf u, \esssup u)$, there exist $s_t\geq 0$, $a_t>0$, and $x_t\in\R^n$  such that the properties (a)--(e) from Theorem \ref{t:klimov-nec-new} are fulfilled.
Then equality holds in the inequality \eqref{ineq_klimov}.
\end{thm}

We conclude this section with two propositions describing quite different situations where equality holds in \eqref{ineq_klimov}. Their proof rests upon
Theorem \ref{t:klimov-suff}.
The first one tells us that, whatever $K$ is, equality holds in \eqref{ineq_klimov} provided that the   super-level sets of $u$ are dilates of any given convex body $L$ containing $0$ in its interior, and the sub-level sets of $\Phi$ are dilates of the polar of $L$.

\begin{proposition}\label{t:klimov-suff1}
 Let $\Phi$ be  as in Theorem \ref{t:klimov} and let $u\in V^{1,\Phi}_{\rm d}(\R^n)$.
Assume that there exist  a  convex body $L$ such that $0\in \inte (L)$,  an increasing function $a: [0, \infty) \to [0, \infty)$, a non-increasing function $b: [0, \infty) \to [0, \infty)$, and $x_0\in\R^n$, such that
\begin{align}
    \label{2024-60}
     \{\Phi\leq s\}=-a(s)L^\circ
\end{align}
for $s>0$,
and
\begin{align}
    \label{2024-61}
     \{u\geq t\}=b(t) L+x_0 \quad \text{ up to a set of Lebesgue measure zero,}
\end{align}
for a.e. $t\in(\essinf u, \esssup u)$.
Then equality holds in the inequality \eqref{ineq_klimov}.
\end{proposition}

The next result shows that the conditions \eqref{2024-60} -- \eqref{2024-61}
 are by no means necessary for equality to hold in \eqref{ineq_klimov}. Indeed,
 it
tells us that
equality may occur
 in \eqref{ineq_klimov} even if none of the level sets of $u$ is homothetic to  another one. For instance, this is the case when $u(x)$ is obtained from a truncation of $-\Phi_\bullet (-x)$. To be more specific, given $t_1 < t_2$, define the truncation function $\cutoff{t_1}{t_2} : \R \to \R$ at the levels $t_1$ and $t_2$ as
$$\cutoff{t_1}{t_2}(t)=
\begin{cases}
t_1, & \text{if $t\leq t_1$,}\\
t, & \text{if $t_1\leq t \leq t_2$,}\\
t_2, & \text{if $t\geq t_2$.}
\end{cases}$$
Then we have what follows.

\begin{proposition}\label{t:klimov-suff2}
Assume that  $\Phi$ is as in Theorem \ref{t:klimov}, fulfills the condition \eqref{2024-6} and
\begin{equation}\label{growth_infinity}
\lim_{|\xi|\to\infty}\frac{\Phi(\xi)}{|\xi|}=\infty.
\end{equation}
Let $u: \rn \to \R$ be any function having the form
\begin{align}
    \label{2024-70prima}
u(x)=\cutoff{t_1}{t_2}\left(t_3-a\,\Phi_\bullet\left(\frac{x_0-x}{a}\right)\right),
\end{align}
for some $t_i\in \R$, $i=1,2,3$, with $t_1 <t_2 \leq t_3$, some $a>0$, and $x_o\in \rn$.
Then $u \in V^{1,\Phi}_{\rm d}(\rn)$ and equality holds in the inequality \eqref{ineq_klimov}.
\end{proposition}

\begin{remark} {\rm The condition~\eqref{growth_infinity} ensures that $\Phi_\bullet$ is finite-valued. The conclusion of Proposition \ref{t:klimov-suff2} continues to hold without \eqref{growth_infinity} provided that  the
function $\Phi_\bullet$  is finite in a neighbourhood of the set $\{\Phi_\bullet\leq (t_3-t_1)/a\}$.  }
\end{remark}

\section{Technical lemmas on $n$-dimensional Young functions}\label{techyoung}

We collect here some properties of $n$-dimensional Young functions, their Young conjugates, and their symmetrals with respect to a convex body which are critical in the proofs of our main results.

We begin with a lemma  on an alternate formula for the Young conjugate of an $n$-dimensional function, via the support functions of its sub-level sets.

\begin{lemma}\label{conj_as_level_sets}
Let $\Phi:\R^n\to [0, \infty]$ be a  lower semicontinuous convex function. Then,
\begin{equation}\label{conjugate_and_support}
   \Phi_\bullet(\xi)=  \sup_{s\geq0}\big(h_{\{\Phi\leq s\}}(\xi)-s\big) \quad \text{ for $\xi \in \rn$.}
 \end{equation}
\end{lemma}
\begin{proof} Fix $\xi \in \rn$.
We claim that
\begin{equation}\label{cls1}
\sup_{\eta\in \R^n}\big(\sprod{\xi}{\eta}-\Phi(\eta)\big) =\sup_{s\geq0, \eta\in\{\Phi\leq s\}}\big(\sprod{\xi}{\eta}-s\big).
\end{equation}
Indeed, on one hand, one has that
\begin{align*}
\sup_{\eta\in \R^n}\big(\sprod{\xi}{\eta}-\Phi(\eta)\big)&=\sup_{s\geq0, \eta\in\{\Phi=s\}}\big(\sprod{\xi}{\eta}-s\big) \leq\sup_{s\geq0, \eta\in\{\Phi\leq s\}}\big(\sprod{\xi}{\eta}-s\big).
\end{align*}
On the other hand, let $\{\eta_k\}$ and $\{s_k\}$ be sequences such that
$\Phi(\eta_k)\leq s_k$ and
$\lim_k \big(\sprod{\xi}{\eta_k}-s_k\big)=\sup_{s\geq0, \eta\in\{\Phi\leq s\}}\big(\sprod{\xi}{\eta}-s\big)$. Then,
\begin{align*}
\sup_{s\geq0, \eta\in\{\Phi\leq s\}}\big(\sprod{\xi}{\eta}-s\big)= \lim_k \big(\sprod{\xi}{\eta_k}-s_k\big) \leq \limsup_k \big(\sprod{\xi}{\eta_k}-\Phi(\eta_k)\big) \leq \sup_{\eta\in \R^n}\big(\sprod{\xi}{\eta}-\Phi(\eta)\big).
\end{align*}
The equality \eqref{cls1} is thus established. Equation \eqref{conjugate_and_support} hence follows, via the following chain:
\begin{align*}
\Phi_\bullet(\xi)&=\sup_{\eta\in \R^n}\big(\sprod{\xi}{\eta}-\Phi(\eta)\big) =\sup_{s\geq0, \eta\in\{\Phi\leq s\}}\big(\sprod{\xi}{\eta}-s\big)  =\sup_{s\geq0}\big(
\sup_{\eta\in\{\Phi\leq s\}}
\sprod{\xi}{\eta}-s\big) =\sup_{s\geq0}\big(h_{\{\Phi\leq s\}}(\xi)-s\big).
\end{align*}
\end{proof}

The next result provides us with an expression for any $n$-dimensional Young whose sub-level sets are homothetic to each other, for its Young conjugate and for the support function of the sub-level sets of the latter.

\begin{lemma}
    \label{convesymm}
    Let $\Phi$ be an $n$-dimensional Young function.  Assume that there exist a convex body $L$ with $0\in \inte(L)$ and a non-decreasing function $a: (0, \infty) \to (0, \infty)$ such that
    \begin{align}
        \label{2024-70}
        \{\Phi \leq s\} = a(s) L^\circ \qquad \text{for $s>0$,}
    \end{align}
    where $L^\circ$ is the polar body of $L$.
    Then, there exists a Young function $A$ such that
\begin{align}
    \label{2024-71}
    \Phi (\xi) = A(H^\smallL(\xi)) \qquad \text{for $\xi \in \rn$,}
\end{align}
where $H^\smallL$ denotes the support function of $L$, and
\begin{align}
    \label{2024-75}
    \Phi_\bullet (\xi) = A_\bullet(H_0^\smallL(\xi)) \qquad \text{for $\xi \in \rn$,}
\end{align}
 where $H_0^\smallL$ is the gauge function of $L$.
\\ Moreover,
\begin{align}
    \label{2024-85bis}
    h_{\{\Phi\leq s\}} (\xi)= A^{-1}(s)H^\smallL_0(\xi) \qquad \text{for $s \geq 0$ and  $\xi \in \rn$}.
\end{align}
and
\begin{align}
    \label{2024-85}
    h_{\{\Phi_{\bullet}\leq s\}} (\xi)= A_\bullet^{-1}(s)H^\smallL(\xi) \qquad \text{for $s \geq 0$ and  $\xi \in \rn$}.
\end{align}

\end{lemma}
\begin{proof}
  Equations \eqref{polarL=hLlq1} and \eqref{2024-70} imply, for $s>0$,
  \begin{align}
      \label{2024-100}
      \{\xi \in \rn: \Phi (\xi)\leq s\} & = a(s)\{\xi \in \rn: H^\smallL(\xi)\leq 1\} \\ \nonumber  & =\{\xi \in \rn: H^\smallL(\xi)\leq a(s)\}= \{\xi \in \rn: a^{-1}(H^\smallL(\xi))\leq s\},
  \end{align}
  where $a^{-1}$ denotes the generalized left-continuous inverse of $a$. Hence,
  \begin{align}
      \label{2024-101}
      \Phi(\xi) = a^{-1}(H^\smallL(\xi)) \qquad \text{for $\xi \in \rn$.}
  \end{align}
  This shows that equation \eqref
{2024-71} holds with $A(t)=a^{-1}(t)$. The fact that the latter function is actually a Young function is a consequence of the assumption that $\Phi$ is an $n$-dimensional Young function.
\\
Equation \eqref{2024-71} implies  \eqref{2024-75},  thanks to \cite[Theorem 15.3]{Rockafellar}.
\\  Equation \eqref{2024-85bis} follows from the chain:
\begin{align}
    \label{2024-104bis}
    h_{\{\Phi\leq s\}} (\xi) & = \sup_{ \{\Phi(\eta)\leq s\}} \sprod{\xi}{\eta} = \sup_{ \{A(H^\smallL(\eta))\leq s\}} \sprod{\xi}{\eta}
      \\ \nonumber & = \sup_{ \{H^\smallL(\eta)\leq 1\}} \sprod{\xi}{A^{-1}(s)\eta}= A^{-1}(s) \sup_{ \{H^\smallL(\eta)\leq 1\}} {\sprod{\xi}\eta}= A^{-1}(s) H^\smallL_0 (\xi)
\end{align}
for $\xi \in \rn$. Here, $A^{-1}$ denotes the generalized right-continuous inverse of $A$. Note that    the last equality relies upon \eqref{L=h0Llq1}. \\ Equation \eqref{2024-85} is a consequence of a parallel chain of equalities.
  \end{proof}

The following lemma informs us on how the non-vanishing or the finiteness of an $n$-dimensional Young function are reflected on the behaviour near zero or near infinity of its Young conjugate.

\begin{lemma}
    \label{limitiPhi}
   Let $\Phi$ be an $n$-dimensional Young function.   Then:
\begin{equation}\label{2024-7'}
\text{if $\Phi(\xi)>0$ for $\xi \neq 0$, \,\, then $\lim_{|\xi|\to 0} \frac{\Phi_\bullet(\xi)}{|\xi|}=0$ \quad \text{and} \quad $\lim_{|\xi|\to 0} \frac{\Phi_{\bullet \smallK}(\xi)}{|\xi|}=0$;}
\end{equation}
\begin{equation}\label{2024-7''}
\text{if $\Phi(\xi)<\infty$ for $\xi \in \rn$, \,\, then
$\lim_{|\xi|\to \infty} \frac{\Phi_\bullet(\xi)}{|\xi|}=\infty$
\quad \text{and} \quad
$\lim_{|\xi|\to \infty} \frac{\Phi_{\bullet \smallK}(\xi)}{|\xi|}=\infty$.}
\end{equation}
\end{lemma}
\begin{proof} We begin by proving the limits concerning the function $\Phi_\bullet$.
    We have that
$$\frac{\Phi_\bullet(\xi)}{|\xi|} = \sup\bigg\{\frac{\sprod{\xi}{\eta}}{|\xi|} - \frac{\Phi(\eta)}{|\xi|}: \eta \in \rn\bigg\}\leq \sup\bigg\{|\eta| - \frac{\Phi(\eta)}{|\xi|}: \eta \in \rn\bigg\}.$$
Let $c, r>0$ be such that $\Phi(\eta)\geq c|\eta|$ if $|\eta|\geq r$. Fix any $\varepsilon \in (0,r)$. Therefore,
\begin{align}
    \label{2024-070}
\frac{\Phi_\bullet(\xi)}{|\xi|} \leq \max\bigg\{\sup\bigg\{|\eta| - \frac{\Phi(\eta)}{|\xi|}: |\eta|< \varepsilon\bigg\},
    \sup\bigg\{|\eta| - \frac{\Phi(\eta)}{|\xi|}: \varepsilon \leq|\eta|\leq  r\bigg\},
    \sup\bigg\{|\eta| - \frac{\Phi(\eta)}{|\xi|}: |\eta|\geq   r\bigg\}
    \bigg\}.
\end{align}
If $|\xi|<c$, then
\begin{align}
    \label{2024-071}
    \sup\bigg\{|\eta| - \frac{\Phi(\eta)}{|\xi|}: |\eta|\geq   r\bigg\}\leq 0.
\end{align}
On the other hand,
\begin{align}
    \label{2024-072}
 \sup\bigg\{|\eta| - \frac{\Phi(\eta)}{|\xi|}: |\eta|< \varepsilon\bigg\}<\varepsilon
\end{align}
for every $\xi \in \rn \setminus \{0\}$.
\\ Finally, since we are assuming that $\Phi (\eta)>0$ for $\eta \neq 0$, we have that
$$\inf_{|\eta|\geq \varepsilon}\frac{\Phi(\eta)}{r-\varepsilon}>0.$$
Thus, if $|\xi| \leq \inf_{|\eta|\geq \varepsilon}{\Phi(\eta)/(r-\varepsilon)}$, then
$$\inf_{r\geq |\eta|\geq \varepsilon}\frac{\Phi(\eta)}{|\eta|-\varepsilon} \geq \inf_{|\eta|\geq \varepsilon}\frac{\Phi(\eta)}{r-\varepsilon}>0,$$
whence
\begin{align}
    \label{2024-073}
    \sup\bigg\{|\eta| - \frac{\Phi(\eta)}{|\xi|}: \varepsilon \leq|\eta|\leq  r\bigg\}<\varepsilon.
\end{align}
Combining equations \eqref{2024-070}--\eqref{2024-073} tells us that
$$\frac{\Phi_\bullet(\xi)}{|\xi|}<\varepsilon$$
provided that $|\xi|< \min\{c, \inf_{|\eta|\geq \varepsilon}\frac{\Phi(\eta)}{r-\varepsilon}\}$. This establishes the first limit in \eqref{2024-7'}.
\\ Next, fix $M>0$ and  set
$$K=\sup_{|\eta|=M}\Phi(\eta).$$
Observe that the supremum is  finite since we are assuming that $\Phi$ is finite-valued. If $|\xi|\geq \frac {2K}M$, then
\begin{align}
    \label{2024-075}
\frac{\Phi_\bullet(\xi)}{|\xi|} = \sup\bigg\{\frac{\sprod{\xi}{\eta}}{|\xi|} - \frac{\Phi(\eta)}{|\xi|}: \eta \in \rn\bigg\} \geq
\sup\bigg\{\frac{\sprod{\xi}{\eta}}{|\xi|} - \frac{\Phi(\eta)}{|\xi|}: |\eta| =M\bigg\}
  \geq M- \frac K{|\xi|} \geq \frac M2,
\end{align}
where the second inequality holds since
$$\frac{\sprod{\xi}{\eta}}{|\xi|} - \frac{\Phi(\eta)}{|\xi|}= M-\frac{\Phi(\eta)}{|\xi|} \quad \quad \text{if $\eta = \frac{\xi}{|\xi|}M$.}$$
The first limit in \eqref{2024-7''} hence follows.
\\ Now, consider the limits involving the function $\Phi_{\bullet \smallK}$. Let $H_0^\smallK$ be the gauge function of the set $K$. Hence, the sub-level sets of $\Phi_{\bullet \smallK}$ are dilates of $\{H_0^\smallK \leq 1\}$. Define the function $\Psi : \rn \to [0, \infty]$ as
$$\Psi (\xi)= \sup_{H_0^\smallK (\eta)\leq H_0^\smallK (\xi)} \Phi_\bullet (\eta) \qquad \text{for $\xi \in \rn$.}$$
Inasmuch as the sub-level sets of $\Psi$ are also dilates of $\{H_0^\smallK \leq 1\}$, we have that $\Psi_\smallK = \Psi$. Moreover, $\Phi_\bullet \leq \Psi$. Therefore, since the operation of symmetrization with respect to $K$ is monotone,
\begin{align}
    \label{2024-080}
    \Phi_{\bullet \smallK}(\xi) \leq \Psi (\xi) \qquad \text{for $\xi \in \rn$.}
\end{align}
Owing to the property \eqref{equivH},
\begin{align}
    \label{2024-081}
    \lim_{|\xi|\to 0} \frac{\Phi_\bullet(\xi)}{|\xi|}=0 \quad \text{if and only if} \quad \lim_{H^\smallK_0(\xi)\to 0} \frac{\Phi_\bullet(\xi)}{H^\smallK_0(\xi)}=0 \quad \text{if and only if} \quad
\lim_{H^\smallK_0(\xi)\to 0} \frac{\Psi(\xi)}{H^\smallK_0(\xi)}=0,
\end{align}
and
\begin{align}
    \label{2024-082}
\lim_{H^\smallK_0(\xi)\to 0} \frac{\Phi_{\bullet \smallK}(\xi)}{H^\smallK_0(\xi)}=0 \quad \text{if and only if} \quad  \lim_{|\xi|\to 0} \frac{\Phi_{\bullet \smallK}(\xi)}{|\xi|}=0.
\end{align}
The second limit in \eqref{2024-7'} follows from the first limit, via \eqref{2024-080} --\eqref{2024-082}.
\\ An analogous argument,
with the replacement of the function $\Psi$ with the function $\Theta$ given by
$$\Theta (\xi)= \inf_{H_0^\smallK (\eta)\geq H_0^\smallK (\xi)} \Phi_\bullet (\eta) \qquad \text{for $\xi \in \rn$,}$$
 shows that the first limit in \eqref{2024-7''} implies the second one.
\end{proof}

Properties of the support function of the sub-level sets of (the Young conjugate of) an $n$-dimensional Young function and of its symmetral with respect to a convex body are established in the last two lemmas of this section.

\begin{lemma}\label{properties_of_Phi}
Let $\Phi$ be an $n$-dimensional Young function.
\\ (i)
 Assume that   $\xi\in\R^n\setminus \{0\}$. Then, the function
\begin{align}
    \label{2024-48}
    [0,+\infty) \ni s\to h_{\{\Phi_{\bullet\smallK}\leq s\}}(\xi)
\end{align}
 is  concave, continuous, nonnegative and nondecreasing.
\\ (ii)   Let $s_0\geq0$ be such that the set  $\{\Phi_\bullet\leq s_0\}$ has  nonempty interior.
Given $\xi\in\R^n\setminus \{0\}$, assume that    $\eta \in \partial \{\Phi_\bullet\leq s_0\}$ is such that $\xi$ belongs to the outer normal cone to $\{\Phi_\bullet\leq s_0\}$ at $\eta$.
Then,  the function
$$[0,+\infty) \ni s\to h_{\{\Phi_{\bullet}\leq s\}}(\xi)-s$$
attains its maximum at $s_0$
if and only if
\[
\Phi_\bullet(\eta)=s_0\quad\text{and}\quad \xi\in\partial\,\Phi_\bullet(\eta).
\]
\end{lemma}
\begin{proof}
Part (i). The function $\Phi_{\bullet \smallK}:\R^n\to[0,+\infty]$ is lower semicontinuous and convex. Consequently, its epigraph
\[
E=\{(\xi,t)\in\R^n\times\R : t\geq\Phi_{\bullet \smallK}(\xi)\}
\]
is closed and convex.
Define the function $\nu : [0, \infty) \to [0, \infty)$ as
\[
\nu(s)=|\{\Phi_{\bullet \smallK}\leq s\}|^{1/n} \qquad \text{for $s \geq 0$.}
\]
Hence,
$$\nu (s)=\cH^n\big(E\cap \{(\xi,t): t=s\}\big)^{1/n} \qquad \text{for $s \geq 0$.}$$
The function $\nu$ is trivially  nondecreasing and nonnegative.
Moreover, given $s_1,s_2\geq0$ and $\lambda\in[0,1]$, the convexity of $E$ guarantees that
\[
\lambda\big(E\cap \{(\xi,t): t=s_1\}\big) +(1-\lambda)\big(E\cap \{(\xi,t): t=s_2\}\big)\subset E\cap \{(\xi,t): t=\lambda s_1+(1-\lambda)s_2\}.
\]
The Brunn-Minkowski inequality for convex bodies \cite[Theorem 7.1]{Sch93}  then implies that
\begin{align*}
\lambda\nu(s_1)+(1-\lambda) \nu(s_2)&\leq
\cH^n\Big( \lambda\big(E\cap \{(\xi,t): t=s_1\}\big) +(1-\lambda)\big(E\cap \{(\xi,t): t=s_2\}\big)\Big)^{1/n}\\
&\leq \nu(\lambda s_1+(1-\lambda)s_2).
\end{align*}
This shows that the function $\nu$ is concave, and hence continuous in $(0, \infty)$. Also, since
 $E$ is closed,  the function $\nu$ is  continuous also at $s=0$.
By the definition of the function $\Phi_{\bullet\smallK}$,
\[
\{\Phi_{\bullet\smallK}\leq s\}=-\nu(s) \kappa^{-1/n} K \qquad \text{for $s \geq 0$,}
\]
where $\kappa=|K|$.
Thus,
\begin{align}
    \label{2024-45}
    h_{\{\Phi_{\bullet\smallK}\leq s\}}(\xi)=\nu(s)\kappa^{-1/n} h_{-K}(\xi) \qquad \text{for $s \geq 0$.}
\end{align}
Notice that $h_{-K}(\xi)>0$,
since the interior of $K$ contains $0$.
Equation \eqref{2024-45}  implies that the function
defined  by \eqref{2024-48} inherits the properties of $\nu$, whence the conclusion follows.
\\ Part (ii).
   The definition of $\eta$ implies  that
\begin{equation}\label{2024-50}
\sprod{\xi}{\eta}=\sup_{\theta \in\{\Phi_{\bullet}\leq s_0\}}\sprod{\xi}{\theta}=h_{\{\Phi_{\bullet}\leq s_0\}}(\xi).
\end{equation}
First, assume that $s_0$ is a maximum point for
$h_{\{\Phi_{\bullet}\leq s\}}(\xi)-s$. The lower semicontinuity of $\Phi_\bullet$ implies $\Phi_\bullet(\eta)\leq s_0$. Therefore,
\begin{equation}\label{ineq_maxima}
    \sprod{\xi}{\eta}-\Phi_\bullet(\eta)\geq h_{\{\Phi_{\bullet}\leq s_0\}}(\xi)-s_0
=\sup_{s\geq0}\big(h_{\{\Phi_\bullet\leq s\}}(\xi)-s\big)\\
    =\Phi(\xi),
\end{equation}
where  the last equality is a consequence of \eqref{conjugate_and_support}.
Thus $\Phi(\xi)+\Phi_\bullet(\eta)\leq\sprod{\xi}{\eta}$, whence, via
\cite[Theorem 23.5]{Rockafellar}, we deduce that $\xi\in\partial\,\Phi_\bullet(\eta)$ and $\Phi(\xi)+\Phi_\bullet(\eta)=\sprod{\xi}{\eta}$.
Coupling the latter equality with equation \eqref{ineq_maxima}
yields $\Phi_\bullet(\eta)=s_0$.
\\ Conversely,
assume  that $\Phi_\bullet(\eta)=s_0$ and $\xi\in\partial\,\Phi_\bullet(\eta)$. \cite[Theorem 23.5]{Rockafellar} implies that $\Phi(\xi)+\Phi_\bullet(\eta)=\sprod{\xi}{\eta}$. Thereby, thanks to \eqref{2024-50} again,
\begin{align*}
    \sup_{s\geq0}\big(h_{\{\Phi_\bullet\leq s\}}(\xi)-s\big)=\Phi(\xi) =\sprod{\xi}{\eta}-\Phi_\bullet(\eta) =h_{\{\Phi_{\bullet}\leq s_0\}}(\xi)-s_0,
\end{align*}
namely, $s_0$ is a maximum point for $h_{\{\Phi_{\bullet}\leq s\}}(\xi)-s$.
\end{proof}
\begin{lemma}
    \label{h}
    Assume that  $\Phi$ is an $n$-dimensional Young function satisfying the assumption \eqref{2024-6}. Given $\xi\in\R^n\setminus \{0\}$,
let  $\varphi : [0, \infty) \to \R$ be function defined  by
\begin{align}
    \label{2024-50'}
    \varphi (s)= h_{\{\Phi_{\bullet\smallK}\leq s\}}(\xi)-s \quad \text{for $s\geq 0$}.
\end{align}
Then, there exists $\max_{s\in [0, \infty)} \varphi (s)$.
\\ In particular, if this maximum  is attained at $0$, then $\varphi(0)>0$,     $\inte (\{\Phi_{\bullet\smallK}\leq 0\})\neq \emptyset$, and $\inte(\{\Phi_{\bullet}\leq 0\})\neq \emptyset$.
\end{lemma}
\begin{proof}
Let $\xi\in\R^n\setminus \{0\}$.  An application of Lemma~\ref{convesymm},  with $K^\circ$ and $\Phi_{\bullet\smallK}$ playing the role of $L$ and $\Phi$, tell us that there exists a Young function $A$ such that
\begin{equation}\label{convesymm_adapted}
\Phi_{\bullet\smallK}(\xi)=A(H_0^\smallK(\xi))
\quad\text{ and }\quad
h_{\{\Phi_{\bullet\smallK}\leq s\}}(\xi)=A^{-1}(s)\,H^\smallK(\xi).
\end{equation}
Here, we have used the standard properties  $H^{\smallK^\circ}=H_0^\smallK$ and $H_0^{\smallK^\circ}=H^\smallK$.
\\
Since $0\in\inte(K)$, we have that $H^\smallK$ is actually a gauge function and hence
 \eqref{equivH} holds. Thanks to \eqref{equivH} and   the first equation in \eqref{convesymm_adapted}, Lemma~\ref{limitiPhi} implies that
\begin{equation*}
 \lim_{t\to0^+}\frac{A(t)}{t}=0\quad \text{ and }\quad    \lim_{t\to+\infty}\frac{A(t)}{t}=\infty.
\end{equation*}
Via the second equation in \eqref{convesymm_adapted}, these limits in turn imply that
\begin{equation*}
    \lim _{s\to 0^+}\frac{h_{\{\Phi_{\bullet\smallK}\leq s\}}(\xi)}s=\infty
\qquad \text{and} \qquad
    \lim _{s\to +\infty}\frac{h_{\{\Phi_{\bullet\smallK}\leq s\}}(\xi)}s=0.
\end{equation*}
As a consequence, the function $\varphi$ defined by \eqref{2024-50'}
enjoys the following properties:
$$\varphi (0)\geq0, \qquad \varphi (s)>0 \,\, \text{for $s$
in a right neighbourhood of $0$}, \qquad   \lim_{s\to +\infty}\varphi (s) = -\infty.$$
By Lemma~~\ref{properties_of_Phi}, Part (i), 
the function $\varphi(s)$ is also concave and continuous in $[0,+\infty)$. Hence, there exists $\max_{s\in [0, \infty)} \varphi (s)$.
\end{proof}

\section{Technical lemmas on Sobolev functions and their symmetrals}

Some key steps of our proof of Theorem \ref{ineq_klimov}, which involve fine properties  of Sobolev functions and of their symmetrals with respect to a convex body, are enucleated in the lemmas of this section.

The first one deals with the derivative of the distribution function  $\mu$ of a Sobolev function $u$,  defined
as in \eqref{mu}.

\begin{lemma}\label{lemma_Andrea}
 Assume that $u\in V^{1,1}_d(\R^n)$.  Then $\Tu\in V^{1,1}_d(\R^n)$ and
\begin{equation}\label{may0}
\int_{\{u=t\}} \frac{d \cH^{n-1}}{|\nabla u|}\,\leq \int_{\{\Tu=t\}} \frac{d \cH^{n-1}}{|\nabla \Tu|} \,  = -\mu '(t) <\infty \quad\quad \text{for a.e. $t>\essinf u$.}
\end{equation}
\end{lemma}
\begin{proof}
The fact that $\Tu\in V^{1,1}_d(\R^n)$ is established in \cite{AFLT}.
As for the inequality \eqref{may0}, observe that
$$\mu (t)= |\{\nabla u =0, u>t\}| + \int_{\{u>t\}} \frac{\chi_{\{\nabla u\neq 0\}}}{|\nabla u|}|\nabla u|\, dx \quad \text{for  $t>\essinf u$.}$$
The function $\mu$ is nonincreasing. Hence, it is differentiable a.e. and, owing to the coarea formula,\eqref{coarea},
\begin{align}\label{may1}
- \mu '(t) & = - \frac{d}{dt} |\{\nabla u =0, u>t\}| - \frac{d}{dt} \int_t^\infty \int_{\{u=\tau\}}\frac{\chi_{\{\nabla u\neq 0\}}(x)}{|\nabla u(x)|} d \cH^{n-1}(x)\, d\tau \\ \nonumber & \geq \int_{\{u=t\}}\frac{\chi_{\{\nabla u\neq 0\}}(x)}{|\nabla u(x)|} d \cH^{n-1}(x)  = \int_{\{u=t\}}\frac{d \cH^{n-1}(x)}{|\nabla u(x)|}
\quad \text{for a.e. $t>\essinf u$.}
\end{align}
  Observe that the last equality holds inasmuch as, for a.e. $t>\essinf u$,
\begin{align}
    \label{2024-41}
    \nabla u(x) \neq 0 \qquad \text{for $\mathcal H^{n-1}$-a.e. $x\in\{u=t\}$.}
\end{align}
Equation \eqref{2024-41} is in turn a consequence of the coarea formula, which implies that
$$0=\int_{\rn}\chi_{\{\nabla u= 0\}}|\nabla u|\, dx = \int_{-\infty}^\infty \int_{\{u=t\}}\chi_{\{\nabla u= 0\}}\, d\mathcal H^{n-1}\, dt= \int_{-\infty}^{\infty} \mathcal H^{n-1}(\{\nabla u= 0\}\cap \{u=t\})\, dt.$$
Denote by $H_0$  the gauge function of $K$. Thereby,
$$K = \{H_0 \leq 1\},$$
and, setting $\kappa = |K|$, we have that
\begin{equation}\label{may6}
 |\{H_0 \leq r\}| =|r K|= k r^n \qquad \text{for $r>0$.}
\end{equation}
Moreover,
\begin{equation}\label{uK_and_gauge}
\Tu (x) = u^*(\kappa H_0(x)^n) \qquad \text{for a.e. $x\in \R^n$.}
\end{equation}
 Since the functions $u$ and $\Tu$ share the same distribution function $\mu$, equation \eqref{may1} also holds with $u$ replaced with $u^\smallK$. Namely,
\begin{align}\label{may2}
- \mu '(t) & = - \frac{d}{dt} |\{\nabla \Tu =0, \Tu>t\}| - \frac{d}{dt} \int_t^\infty \int_{\{\Tu=\tau\}}\frac{\chi_{\{\nabla \Tu\neq 0\}}(x)}{|\nabla \Tu(x)|} d \cH^{n-1}(x)\, d\tau \\ \nonumber & = - \frac{d}{dt} |\{\nabla \Tu =0, \Tu>t\}| +  \int_{\{\Tu=t\}}\frac{d \cH^{n-1}(x) }{|\nabla \Tu(x)|} \quad \text{for a.e. $t>\essinf u$.}
\end{align}
The inequality \eqref{may0} will follow from \eqref{may1} and \eqref{may2} if we show that
\begin{equation}\label{may3}
\frac{d}{dt} |\{\nabla \Tu =0, \Tu>t\}| =0\quad\quad \text{for a.e. $t>\essinf u$.}
\end{equation}
To prove  equation \eqref{may3}, let us begin by observing that, if $E$ is a measurable subset of $[0, \infty)$, then
\begin{equation}\label{may4}
|E| = |\{x \in  \R^n: \kappa H_0(x)^n \in E\}|.
\end{equation}
Indeed,
\begin{align}\label{may5}
 |\{x \in  \R^n: \kappa H_0(x)^n \in E\}
 & = \int_{\R^n} \chi_E(\kappa H_0(x)^n)\, dx=  \int_{\R^n} \chi_E\left(|x|^n\kappa H_0\left(\frac x{|x|}\right)^n\right)\, dx\\ \nonumber
& = \int_0^\infty \int_{\mathbb S^{n-1}} \chi_E\big(r^n\kappa H_0(\nu)^n\big)d \cH^{n-1}(\nu) r^{n-1}\, dr \\ \nonumber
&=  \frac1{n k}\int_0^\infty \int_{\mathbb S^{n-1}} \chi_E(s)\frac{d \cH^{n-1}(\nu)}{H_0(\nu)^n} \, ds \\ \nonumber
& =\frac{|E|}{n\kappa} \int_{\mathbb S^{n-1}} \frac{d \cH^{n-1}(\nu)}{H_0(\nu)^n}.
\end{align}
An application of equation \eqref{may5} with $E=[0, 1]$ and equation \eqref{may6} tell us that
\begin{align*}
1=\frac{1}{n\kappa} \int_{\mathbb S^{n-1}} \frac{d \cH^{n-1}(\nu)}{H_0(\nu)^n}.
\end{align*}
Equation \eqref{may4} follows from the latter equality and \eqref{may5}.
\\  { Since $\Tu$ is the composition of a convex function with a nonincreasing locally absolutely continuous function, we have that}
\begin{equation}\label{giu1}
\nabla \Tu (x)= (u^*)'(\kappa H_0^n(x)) \, n\kappa H_0(x)^{n-1}\nabla H_0(x)
\end{equation}
for every $x\in \R^n$ such that $(u^*)'(\kappa H_0^n(x))$ and $\nabla H_0(x)$ exist, and hence
for a.e. $x\in \R^n$. Moreover,
$$ H_0(x)^{n-1}\nabla H_0(x)\neq 0 \qquad
\text{ for a.e $x\in \rn$. }
$$
Therefore,
$$|\{x\in \R^n: \nabla \Tu (x) =0, \Tu (x) >t\}|= |\{x\in \R^n: (u^*)'(\kappa H_0^n(x)) =0,  u^*(\kappa H_0^n(x)) >t\}| \quad \text{for $t > \essinf u$}.$$
Hence, on choosing
$$E =\{s>0:  (u^*)'(s)=0, u^*(s)>t\}$$
from equation \eqref{may4} we infer  that
\begin{equation}\label{may7}
\frac{d}{dt} |\{\nabla \Tu =0, \Tu>t\}| = \frac{d}{dt}|\{(u^*)'(s)=0, u^*(s)>t\}| \quad \text{for a.e. $t>\essinf u$.}
\end{equation}
 On the other hand, \cite[Lemma 2.4]{CianchiFuscoARMA} entails that
\begin{equation}\label{may8}
 \frac{d}{dt}|\{(u^*)'(s)=0, u^*(s)>t\}|=0 \quad \text{for a.e. $t>\essinf u$.}
\end{equation}
Equation \eqref{may3} is a consequence of \eqref{may7} and \eqref{may8}.
\end{proof}

Loosely speaking, the next lemma amount to an application of the isoperimetric inequality \eqref{may9} to the case when $E$ is a super-level set of  a Sobolev function $u$. A delicate point in this application is the identification of the unit normal vector $\nu^{E}$ with $-\frac{\nabla u}{|\nabla u|}$. This requires sophisticated results from Geometric Measure Theory.

\begin{lemma}\label{lemma_Andrea1}
 Let $L\subset\R^n$ be a convex body such that $0\in \inte (L)$ and let $u\in V^{1,1}_d(\R^n)$.  Then
\begin{equation}\label{may10}
\int_{\{u=t\}} h_L\bigg(-\frac{\nabla u}{|\nabla u|}\bigg)\, d \cH^{n-1} \geq n |\{u>t\}|^{\frac{n-1}{n}} |L|^{\frac 1n} \quad\quad \text{for a.e. $t>\essinf u$.}
\end{equation}
Moreover,  for a.e. $t\in (\essinf u, \esssup u)$,  equality holds in \eqref{may10} if and only if
\begin{equation}\label{may11}
\{u>t\} = \{u>t\}^\smallL \qquad \text{up to a translation and up to sets of Lebesgue measure zero.}
\end{equation}
\end{lemma}
\begin{proof}
Thanks to the property \eqref{2024-6prima},
the conclusion will follow from the  isoperimetric inequality \eqref{may9}, 
provided that we show that, for a.e. $t\in (\essinf u, \esssup u)$,
\begin{equation}\label{may12}
\nu^t(x) = -\frac{\nabla u(x)}{|\nabla u(x)|}  \qquad \text{for $\cH^{n-1}$-a.e.   $x\in \partial ^* \{u>t\}$},
\end{equation}
where $\nu^t$ denotes  the geometric measure theoretical outer unit normal to the set $\{u>t\}$.
\\ To verify equation \eqref{may12}, recall that the subgraph $S$ of $u$, defined as
$$S = \{(x,t)\in \mathbb R^{n+1}: x\in \R^n,  u(x)>t\},$$
is a set of finite perimeter in $\R^n$. Moreover, the  outer unit normal $\nu^S$ to $S$ satisfies
\begin{equation}\label{may13}
\nu^S(x,t)= \bigg(- \frac{\nabla u(x)}{\sqrt{1+|\nabla u(x)|^2}}\,, 1\bigg) \quad \text{for $\mathcal H^n$-a.e. $(x,t) \in \partial ^*S$,}
\end{equation}
see e.g. \cite[Theorem C]{CianchiFuscoAdvMath}. On setting
$$S_t = \{x\in \R^n: (x,t) \in S\} \qquad \text{for $t\in (\essinf u, \esssup u)$,}$$
one trivially has that
\begin{align}\label{may15}
S_t = \{u>t\} \qquad \text{for $t\in (\essinf u, \esssup u)$.}
\end{align}
Denote by $\nu^S_x\in \R^n$ the vector of the first $n$ components of $\nu^S$ and by $\nu^{S_t}$ the outer unit normal vector to $S_t$.
Thanks to \cite[Theorem 2.4]{BarchiesiCagnettiFusco},
$$\partial^* S_t = \big(\partial ^*S\big)_t  \qquad \text{up to sets of of $\cH^{n-1}$ measure zero,}$$
for a.e. $t\in (\essinf u, \esssup u)$, and
\begin{align}\label{may14}
\nu^{t}(x) = \frac{\nu^S_x (x,t)}{|\nu^S_x(x,t)|} = -\frac{\nabla u(x)}{|\nabla u(x)|}\qquad \text{for  $\cH^{n-1}$-a.e. $x \in\partial^* S_t$,}
\end{align}
for a.e. $t\in (\essinf u, \esssup u)$.
  Coupling \eqref{may13} with \eqref{may14} implies \eqref{may12}.
\end{proof}

We conclude this section by showing that the support  function of any sub-level set of the function $\Phi_{\bullet\smallK}$ attains a constant value on $\nabla u^\smallK (x)$ when $x$ ranges on a given level set of  $u^\smallK$. As above, here $\Phi$ denotes any $n$-dimensional Young function, $K$ a convex body containing $0$ in its interior, and $u$ a Sobolev function.

\begin{lemma}\label{l:s(t)_well_posed} Let $\Phi$ be an $n$-dimensional Young function.
 Assume that
 $u\in V^{1,1}_d(\R^n)$. Then,
 for a.e.  $t\in (\essinf u, \esssup u)$,
 \begin{align}\label{2024-10}
\text{$\Tu$ is differentiable $\mathcal H^{n-1}$-a.e. in $\{\Tu =t\}$, \quad $\nabla \Tu (x)\neq 0$  for $\mathcal H^{n-1}$-a.e. $x\in \{\Tu =t\}$,}
 \end{align}
and
\begin{equation}\label{same_argmax}
  h_{\{\Phi_{\bullet\smallK}\leq s\}}(\nabla \Tu (x_1))=h_{\{\Phi_{\bullet\smallK}\leq s\}}(\nabla \Tu (x_2))
 \end{equation}
   for every $s\geq0$ and for $\mathcal H^{n-1}$-a.e. $x_1,x_2\in\{\Tu=t\}$.

\end{lemma}
\begin{proof}
Let $u^*$ be the decreasing rearrangement of $u$.
 Define the sets
$$D_{u^*}^+= \{s>0: (u^*)'(s) \,\, \text{exists, and}\,\, (u^*)'(s)\neq 0\} $$
and
$$D_{u^*}^0= \{s>0: (u^*)'(s) \,\, \text{exists, and}\,\, (u^*)'(s)=0\}.$$
Since $u^*$ is differentiable a.e. in $(0,\infty)$, we have that
\begin{equation}\label{2024-1}
|(0, \infty)\setminus (D_{u^*}^+\cup D_{u^*}^0)|=0.
\end{equation}
Plainly,
$$u^*((0,\infty))= u^*(D_{u^*}^+) \cup u^*(D_{u^*}^0)\cup u^*\big((0, \infty)\setminus (D_{u^*}^+\cup D_{u^*}^0)\big).$$
Since the function $u^*$ is absolutely continuous, the property \eqref{2024-1} ensures that
$$\big|u^*\big((0, \infty)\setminus (D_{u^*}^+\cup D_{u^*}^0)\big)\big|=0.$$
From \cite[Equation (2.22)]{CianchiFuscoARMA} , we infer that
$$|u^*(D_{u^*}^0)|=0.$$
Finally, \cite[Equation (3.12)]{CianchiFuscoARMA}  implies that, if $t\in u^*(D_{u^*}^+)$, then $(u^*)^{-1}(t)$ is a singleton, and hence $(u^*)^{-1}(t) \in D_{u^*}^+$. Altogether, we have that
\begin{align}\label{2024-2}
(u^*)^{-1}(t) \in D_{u^*}^+ \qquad \text{for a.e. $t\in (\essinf u, \esssup u)$.}
\end{align}
Next, let $H_0$ denote  gauge function of $K$, and set $\kappa=|K|$. Then,
for every $r>0$,
\begin{align}\label{2024-3}
\text{there exists $\nabla H_0(x)\neq 0$  for $\mathcal H^{n-1}$-a.e. $x \in \{\kappa H_0^n(x)=r\}$. }
\end{align}
From equations \eqref{giu1}, \eqref{2024-2}, and \eqref{2024-3} we deduce that, for a.e.  $t\in (\essinf u, \esssup u)$,
\begin{align}
    \label{2024-33}
    \text{there exists} \,\,\nabla \Tu (x)= (u^*)'(\kappa H_0^n(x)) \, n\kappa H_0(x)^{n-1}\nabla H_0(x)\neq 0 \quad \text{for $\mathcal H^{n-1}$-a.e. $x\in \{\Tu =t\}$.}
\end{align}
The assertion \eqref{2024-10} is thus established.
\\
In view of the properties discussed above, it suffices to prove \eqref{same_argmax} for every $s \geq 0$ and $t \in u^*(D_{u^*}^+)$. Fix any  $$t \in u^*(D_{u^*}^+).$$
 Let  $s\geq 0$. If $s=0$ and  $\{\Phi_{\bullet\smallK}\leq 0\}=\{0\}$, then \eqref{same_argmax} trivially holds, since  both sides vanish. Assume now that $\{\Phi_{\bullet\smallK}\leq s\}\neq \{0\}$.
Let $t\in (\essinf u, \esssup u)$ be such that equation \eqref{2024-10} holds.
 Let  $x_1, x_2 \in  \{\Tu =t\}$ be such that  $\nabla \Tu (x_1)$ and $\nabla \Tu (x_2)$ exist and are different from $0$.  Owing to equation \eqref{2024-33}, we have that
 $(u^*)^\prime(\kappa H_0(x_1)^n)= (u^*)^\prime(\kappa H_0(x_2)^n)\neq0$, and since $H_0(x_1)=H_0(x_2)$,
\[
\nabla \Tu(x_i)=-a\nabla H_0(x_i) \quad \text{for $i=1,2$,}
\]
for some  $a>0$.
On the other hand,
\[
 \{\Phi_{\bullet\smallK}\leq s\}=-b K,
\]
for some $b>0$ and, by the homogeneity of $h_K$ and equation \eqref{nablaH},
$$ h_{\{\Phi_{\bullet\smallK}\leq s\}}(\nabla\Tu(x_i))=h_{-b K}(-a\nabla H_0(x_i))= a b h_K(\nabla H_0(x_i)) =  a b H(\nabla H_0(x_i))=a b \quad \text{for $i=1,2$.}$$
Hence, equation \eqref{same_argmax} follows.
\end{proof}

\section{Proofs of the main results}

Having the technical material established in the preceding sections at disposal, we are now in a position to accomplish the proof of Theorem \ref{t:klimov}.

\begin{proof}[Proof of Theorem~\ref{t:klimov}]
Assume, for the time being, that $u\in V^{1,1}_{\rm d}(\rn)$. Hence, by Lemma \ref{lemma_Andrea},  $u^\smallK\in V^{1,1}_{\rm d}(\rn)$ as well.
We shall prove that
\begin{equation}\label{ineq_klimov_level_lines}
\int_{\{u=t\}}\frac{\Phi(\nabla u)}{|\nabla u|}\,d \cH^{n-1}\geq
\int_{\{\Tu=t\}}\frac{\Phi_{\bullet\smallK \bullet}(\nabla \Tu)}{|\nabla \Tu|}\,d \cH^{n-1} \quad \text{for a.e. $t\in(\essinf u,\esssup u)$.}
\end{equation}
Fix $\varepsilon >0$.
Let $t\in (\essinf u,\esssup u)$ be such that the properties \eqref{2024-10} and \eqref{same_argmax} hold.
First, suppose that
\begin{equation}\label{2024-11}
\sup_{s\geq 0} \big(h_{\{\Phi_{\bullet\smallK}\leq s\}}(\nabla \Tu (x))-s\big)<\infty
\end{equation}
for $\mathcal H^{n-1}$-a.e. $x\in\{\Tu=t\}$. Notice that, by \eqref{same_argmax}, the supremum in \eqref{2024-11} is independent of $x\in\{\Tu=t\}$, up to subsets of $\mathcal H^{n-1}$ measure zero.
\\
Let $s(t,\varepsilon)\in [0, \infty)$ be such that
$$
h_{\{\Phi_{\bullet\smallK}\leq s(t,\varepsilon)\}}(\nabla \Tu (x))-s(t, \varepsilon) \geq (1-\varepsilon)\sup_{s\geq 0} \big(h_{\{\Phi_{\bullet\smallK}\leq s\}}(\nabla \Tu (x))-s\big).$$
Thanks to the formula~\eqref{conjugate_and_support},
\begin{align} \label{step3}
\int_{\{u=t\}}\frac{\Phi (\nabla u)}{|\nabla u|}\,d \cH^{n-1}
&=\int_{\{u=t\}}\frac{\sup_{s\geq0}\big(h_{\{\Phi_\bullet\leq s\}}(\nabla u)-s\big)}{|\nabla u|}\,d \cH^{n-1} \\
 \nonumber &\geq \int_{\{u=t\}}\frac{h_{\{\Phi_\bullet\leq s(t, \varepsilon)\}}(\nabla u)-s(t,\varepsilon)}{|\nabla u|}\,d \cH^{n-1} \\
\nonumber &=\int_{\{u=t\}}h_{-\{\Phi_\bullet\leq s(t, \varepsilon)\}}\left( -\frac{\nabla u}{|\nabla u|}\right)\,d \cH^{n-1}
-\int_{\{u=t\}}\frac{s(t, \varepsilon)}{|\nabla u|}\,d \cH^{n-1}.
\end{align}
From Lemma~\ref{lemma_Andrea1}, the equimeasurability of the rearrangement, and Lemma~\ref{lemma_Andrea1} again, combined with the fact that $\{\Tu\geq t\}$ and $-\{\Phi_{\bullet\smallK}\leq s(t)\}$ are homothetic, one can deduce that
\begin{align}\label{step4}
\int_{\{u=t\}}
h_{-\{\Phi_{\bullet}\leq s(t,\varepsilon)\}}\left( -\frac{\nabla u}{|\nabla u|} \right)\,d \cH^{n-1}&\geq n\, |\{u\geq t\}|^\frac{n-1}{n}\, |\{\Phi_\bullet\leq s(t,\varepsilon)\}|^\frac1{n}\\ \nonumber
&=n\, |\{\Tu\geq t\}|^\frac{n-1}{n}\, |\{\Phi_{\bullet\smallK}\leq s(t,\varepsilon)\}|^\frac1{n}\\
 \nonumber &=\int_{\{\Tu=t\}} h_{-\{\Phi_{\bullet\smallK}\leq s(t,\varepsilon)\}}\left( -\frac{\nabla \Tu}{|\nabla \Tu|} \right)\,d \cH^{n-1}.
\end{align}
On the other hand,
Lemma~\ref{lemma_Andrea} ensures that
\begin{equation}\label{step5}
 -\int_{\{u=t\}}\frac{s(t,\varepsilon)}{|\nabla u|}\,d \cH^{n-1}\geq
-\int_{\{\Tu=t\}} \frac{s(t,\varepsilon)}{|\nabla \Tu|}\,d \cH^{n-1}.
\end{equation}
Finally, by~\eqref{conjugate_and_support} and the definition of $s(t, \varepsilon)$,
\begin{align}\label{step6}
 \int_{\{\Tu=t\}} h_{-\{\Phi_{\bullet\smallK}\leq s(t,\varepsilon)\}}\left( \frac{-\nabla \Tu}{|\nabla \Tu|} \right)\,d \cH^{n-1}
 & -\int_{\{\Tu=t\}} \frac{s(t, \varepsilon)}{|\nabla \Tu|}\,d \cH^{n-1}
 \\ \nonumber & =
\int_{\{\Tu=t\}}\frac{h_{\{\Phi_{\bullet\smallK}\leq s(t, \varepsilon)\}}(\nabla \Tu)-s(t,\varepsilon)}{|\nabla \Tu|}\,d \cH^{n-1}
\\ \nonumber & \geq (1-\varepsilon)
\int_{\{\Tu=t\}}\frac{\sup_{s\geq 0} \big(h_{\{\Phi_{\bullet\smallK}\leq s\}}(\nabla \Tu)-s\big)}{|\nabla \Tu|}\,d \cH^{n-1}
\\ \nonumber &
= (1-\varepsilon)\int_{\{\Tu=t\}}\frac{\Phi_{\bullet\smallK\bullet}(\nabla \Tu)}{|\nabla \Tu|}\,d \cH^{n-1}.
\end{align}
Combining equations \eqref{step3}--\eqref{step6} yields
\begin{equation}\label{ineq_klimov_level_lines-ep}
\int_{\{u=t\}}\frac{\Phi (\nabla u)}{|\nabla u|}\,d \cH^{n-1}\geq (1-\varepsilon)
\int_{\{\Tu=t\}}\frac{\Phi_{\bullet\smallK\bullet}(\nabla \Tu)}{|\nabla \Tu|}\,d \cH^{n-1}.
\end{equation}
The inequality
\eqref{ineq_klimov_level_lines} hence follows, via the arbitrariness of $\varepsilon$.
 \\ In the case when
 $$\sup_{s\geq 0} \big(h_{\{\Phi_{\bullet\smallK}\leq s\}}(\nabla \Tu (x))-s\big)=\infty,$$
 an obvious variant of the argument above tells us that both sides of the inequality \eqref{ineq_klimov_level_lines} equal $\infty$, whence the relevant inequality holds trivially also in this case.
 \\
An integration of \eqref{ineq_klimov_level_lines} with respect to $t$ over $(\essinf u,\esssup u)$ gives  the inequality
\eqref{ineq_klimov}, via the coarea formula.

To remove the temporary assumption $u\in V^{1,1}_{\rm d}(\rn)$, define for $t>\essinf u$ the function $T_t(u)$ as $T_t(u)= \max\{u, t\}$ and recall from \eqref{2024-31} that $T_t(u)\in V^{1,1}_{\rm d}(\rn)$. An application of \eqref{ineq_klimov} with $u$ replaced with $T_t(u)$ yields
\begin{align}\label{ineq_klimov-1}
\int_{\{u^\smallK >t\}}\Phi_{\bullet\smallK \bullet}(\nabla u^\smallK)\, dx & =\int_{\R^n}\Phi_{\bullet\smallK \bullet}(\nabla T_t(u^\smallK))\, dx
=\int_{\R^n}\Phi_{\bullet\smallK \bullet}(\nabla (T_t(u))^\smallK)\, dx
\\ \nonumber & \leq \int_{\R^n} {\Phi}(\nabla T_t(u))\, dx= \int_{\{u>t\}} {\Phi}(\nabla u)\, dx
\end{align}
 for $t>\essinf u$,
where the second equality holds since $T_t(u^\smallK)=(T_t(u))^\smallK$ a.e. in $\rn$. The inequality \eqref{ineq_klimov} follows from \eqref{ineq_klimov-1}, by passing to the limit as $t\downarrow\essinf u$, owing to the monotone convergence theorem.
\end{proof}

The proof of Theorem~\ref{t:klimov-nec-new} basically consists of decoding the information contained in the fact that, if equality holds in \eqref{ineq_klimov}, then each intermediate inequality in the chain which yields \eqref{ineq_klimov} must also hold as an equality. Conversely, the proof of Theorem \ref{t:klimov-suff} requires checking that its assumptions ensure that all the inequalities in the relevant chain hold as equalities, whence the equality in \eqref{ineq_klimov} follows.

\begin{proof}[Proof of Theorem~\ref{t:klimov-nec-new}]
\emph{Part (i)}.
By a truncation argument as in the proof of Theorem \ref{t:klimov}, we may assume, without loss of generality, that $u\in V^{1,1}_{\rm d}(\rn)$. Indeed,
if
\begin{align}
    \label{2024-204}
\int_{\rn}\Phi_{\bullet\smallK \bullet}(\nabla u^\smallK)\, dx   = \int_{\rn} {\Phi}(\nabla u)\, dx,
\end{align}
 then
\begin{align}
    \label{2024-201}
\int_{\R^n}\Phi_{\bullet\smallK \bullet}(\nabla (T_t(u))^\smallK)\, dx
= \int_{\R^n} {\Phi}(\nabla T_t(u))\, dx
\end{align}
as well.
To verify this claim,  recall from \eqref{ineq_klimov-1} that
\begin{align}
    \label{2024-206}
    \int_{\R^n}\Phi_{\bullet\smallK \bullet}(\nabla (T_t(u))^\smallK)\, dx
 \leq \int_{\R^n} {\Phi}(\nabla T_t(u))\, dx,
\end{align}
and
observe that, on setting $v_t= u- T_t(u)$,  one analogously has:
\begin{align}\label{ineq_klimov-2}
\int_{\R^n}\Phi_{\bullet\smallK \bullet}(\nabla (v_t)^\smallK)\, dx
 \leq \int_{\R^n} {\Phi}(\nabla v_t)\, dx.
\end{align}
On the other hand
\begin{align}\label{2024-202}
 \int_{\R^n} {\Phi}(\nabla u)\, dx  = \int_{\R^n} {\Phi}(\nabla T_t(u))\, dx + \int_{\R^n} {\Phi}(\nabla v_t)\, dx
\end{align}
and
\begin{align}\label{2024-203}
\int_{\R^n}\Phi_{\bullet\smallK \bullet}(\nabla u^\smallK)\, dx=
\int_{\R^n}\Phi_{\bullet\smallK \bullet}(\nabla   (T_t(u))^\smallK)\, dx +\int_{\R^n}\Phi_{\bullet\smallK \bullet}(\nabla (v_t)^\smallK)\, dx
\end{align}
Equation \eqref{2024-201} follows via  \eqref{2024-204}
and \eqref{2024-206} -- \eqref{2024-203}.
\\ Now, let $x\in\{\Tu=t\}$. Owing to  \eqref{2024-10} and Lemma \ref{h},
 for a.e. $t\in (\essinf u, \esssup u)$ the supremum in \eqref{2024-11} is attained at some $s(t) \in [0, \infty)$, and
 \begin{equation}
     \label{interior}
     \text{$\inte (\{\Phi_{\bullet \smallK}\leq s(t)\}) \neq \emptyset$.}
 \end{equation}
Consequently,  for a.e. $t\in (\essinf u, \esssup u)$, the formulas \eqref{step3}--\eqref{step6}
in the proof of Theorem \ref{t:klimov} are fulfilled  with $\varepsilon =0$. Namely,  for a.e. $t\in (\essinf u, \esssup u)$ the following chain holds:
\begin{align}
    \label{2024-12}
    \int_{\{u=t\}}\frac{\Phi (\nabla u)}{|\nabla u|}\,d \cH^{n-1}
&=\int_{\{u=t\}}\frac{\sup_{s\geq0}\big(h_{\{\Phi_\bullet\leq s\}}(\nabla u)-s\big)}{|\nabla u|}\,d \cH^{n-1}
\\ \nonumber
&\geq
\int_{\{u=t\}}\frac{h_{\{\Phi_\bullet\leq s(t)\}}(\nabla u)-s(t)}{|\nabla u|}\,d \cH^{n-1} \\
\nonumber &=\int_{\{u=t\}}h_{-\{\Phi_\bullet\leq s(t)\}}\left( -\frac{\nabla u}{|\nabla u|}\right)\,d \cH^{n-1}
-\int_{\{u=t\}}\frac{s(t)}{|\nabla u|}\,d \cH^{n-1}
\\ \nonumber & \geq
n\, |\{u\geq t\}|^\frac{n-1}{n}\, |\{\Phi_\bullet\leq s(t)\}|^\frac1{n} -\int_{\{\Tu=t\}} \frac{s(t)}{|\nabla \Tu|}\,d \cH^{n-1}\\ \nonumber
&=n\, |\{\Tu\geq t\}|^\frac{n-1}{n}\, |\{\Phi_{\bullet\smallK}\leq s(t)\}|^\frac1{n}  - \int_{\{\Tu=t\}} \frac{s(t)}{|\nabla \Tu|}\,d \cH^{n-1}
\\ \nonumber & =
 \int_{\{\Tu=t\}} h_{-\{\Phi_{\bullet\smallK}\leq s(t)\}}\left( -\frac{\nabla \Tu}{|\nabla \Tu|} \right)\,d \cH^{n-1} - \int_{\{\Tu=t\}} \frac{s(t)}{|\nabla \Tu|}\,d \cH^{n-1}
\\ \nonumber & = \int_{\{\Tu=t\}}\frac{h_{\{\Phi_{\bullet\smallK}\leq s(t)\}}(\nabla \Tu)-s(t)}{|\nabla \Tu|}\,d \cH^{n-1}
\\ \nonumber & =
\int_{\{\Tu=t\}}\frac{\max_{s\geq 0} \big(h_{\{\Phi_{\bullet\smallK}\leq s\}}(\nabla \Tu)-s\big)}{|\nabla \Tu|}\,d \cH^{n-1}
\\ \nonumber &
=  \int_{\{\Tu=t\}}\frac{\Phi_{\bullet\smallK\bullet}(\nabla \Tu)}{|\nabla \Tu|}\,d \cH^{n-1}.
\end{align}
Under the assumption that   equality holds in \eqref{ineq_klimov},  the inequalities in \eqref{2024-12} must hold as equalities. In particular, since   equality holds in the second one,  by Lemma \ref{lemma_Andrea1},
for a.e. $t\in (\essinf u, \esssup u)$
equation \eqref{t:klimov-suff-a}  holds with $s_t=s(t)$, for some $a_t>0$ and $x_t \in \mathbb R^n$. Moreover, equation \eqref{2024-30} is fulfilled with the same choice of $s_t$, thanks to \eqref{interior}.
\\ Equation
\eqref{t:klimov-suff-b} follows
via  Lemma \ref{lemma_Andrea}.
Finally, the conditions~\eqref{t:klimov-suff-c} and \eqref{t:klimov-suff-d} follow from Lemma~\ref{properties_of_Phi}, Part (ii), with $\xi=\nabla u(x)$ and   $\xi=\nabla\Tu(x)$, respectively,  since $s(t)$  maximizes both
\[
h_{\{\Phi_{\bullet}\leq s\}}(\nabla u(x))-s\quad\text{and}\quad h_{\{\Phi_{\bullet\smallK}\leq s\}}(\nabla\Tu(x))-s.
\]
It remains to prove that the function $u$ is quasi-convex, namely, that the set $\{u \geq t_0\}$  is convex for every $t_0 \in (\essinf u, \esssup u)$.
 To this purpose,  choose a sequence $\{t_k\}$ such that equation \eqref{t:klimov-suff-a} holds with $t=t_k$ and $t_k \uparrow t_0$. One has that
\begin{align}
    \label{2024-82}
\{u \geq  t_0\} = \bigcap_{k \in \N}\{u \geq t_k\}.
\end{align}
Since the set $\{u \geq t_k\}$ is convex for every $k\in \N$, the set $\{u \geq t_0\}$ is convex as well.
\\ \emph{Part (ii)}.  Since $\Phi$ is  finite-valued, the assumption \eqref{strict} implies  that $\Phi_\bullet$ is differentiable in  $\inte(\dom\Phi_\bullet)$, and that $\partial \Phi_\bullet(\eta)=\{\nabla\Phi_\bullet(\eta)\}$ for $\eta\in\inte(\dom\Phi_\bullet)$, while $\partial \Phi_\bullet(\eta)=\emptyset$ for $\eta\notin\inte(\dom\Phi_\bullet)$ \cite[Theorems 26.1, 26.3]{Rockafellar}. Thus \eqref{t:klimov-suff-c} implies that $\xi\in\inte(\dom\Phi_\bullet)$ and yields \eqref{equality_necessity_grad1}.
\\
If \eqref{diff} also holds, then $\Phi_\bullet$ is strictly convex on $\inte(\dom\Phi_\bullet)$ \cite[Theorem 26.3]{Rockafellar}. This implies that  there exists exactly one point  $\xi$ such that $\nabla \Phi_\bullet (\xi)$ agrees with the prescribed vector $\nabla u(x)$.
\end{proof}

\begin{proof}[Proof of Theorem~\ref{t:klimov-suff}]   As in the previous proofs, we may assume, without loss of generality, that $u\in V^{1,1}_{\rm d}(\rn)$.
Let $t\in(\essinf u, \esssup u)$ be such that the  properties \eqref{t:klimov-suff-a} -- \eqref{t:klimov-suff-b} hold. The properties ~\eqref{t:klimov-suff-c} and~\eqref{t:klimov-suff-d} imply, via Lemma~\ref{properties_of_Phi}, Part (ii), that $s_t$ maximizes both
\[
h_{\{\Phi_{\bullet}\leq s\}}(\nabla u(x))-s\quad\text{and}\quad h_{\{\Phi_{\bullet\smallK}\leq s\}}(\nabla\Tu(y))-s,
\]
for a.e. $x\in\{u=t\}$ and $y\in\{\Tu=t\}$.
Coupling this piece of information  with ~\eqref{t:klimov-suff-a} and ~\eqref{t:klimov-suff-b} implies that all the inequalities in \eqref{2024-12}, with $s(t)=s_t$, hold as equalities. Altogether,
\[
\int_{\{u=t\}}\frac{\Phi(\nabla u)}{|\nabla u|}\,d \cH^{n-1}=
\int_{\{\Tu=t\}}\frac{\Phi_{\bullet\smallK \bullet}(\nabla \Tu)}{|\nabla \Tu|}\,d \cH^{n-1} \quad \text{for a.e. $t\in (\essinf u, \esssup u)$.}
\]
Thanks to the coarea formula, integrating the latter identity over $(\essinf u, \esssup u)$ yields
$$ \int_{\R^n} {\Phi}(\nabla u)\, dx=\int_{\R^n}\Phi_{\bullet \smallK \bullet}(\nabla \Tu)\,
dx.$$
\end{proof}

We conclude with the proofs of Propositions \ref{t:klimov-suff1}  and \ref{t:klimov-suff2}. As mentioned above, they make use of Theorem \ref{t:klimov-suff}.

\begin{proof}[Proof of Proposition~\ref{t:klimov-suff1}]
Since the two sides of the inequality \eqref{ineq_klimov} are translation invariant, we may assume, without loss of generality, that $x_0=0$.
Theorem~\ref{t:klimov} tells us that
\begin{equation}\label{ineq_klimov_r}
\int_{\R^n}\Phi_{\bullet \smallK \bullet}(\nabla \Tu)\, dx\leq \int_{\R^n} {\Phi}(\nabla u)\, dx.
\end{equation}
The same theorem, applied with $K$, $u$ and $\Phi$ replaced by $L$, $\Tu$ and $\Phi_{\bullet\smallK \bullet}$, respectively, implies that
\begin{equation}\label{ineq_klimov_reverse}
\int_{\R^n}\left(\Phi_{\bullet\smallK\bullet}\right)_{\bullet\smallL \bullet}\left(\nabla \left(\Tu\right)^\smallL\right)\, dx\leq \int_{\R^n}\Phi_{\bullet\smallK\bullet}
(\nabla \Tu)\, dx.
\end{equation}
 We claim  that
 \begin{align}
     \label{2024-63}\int_{\R^n}\left(\Phi_{\bullet\smallK\bullet}\right)_{\bullet\smallL \bullet}\left(\nabla \left(\Tu\right)^\smallL\right)\, dx=\int_{\R^n} {\Phi}(\nabla u)\, dx.
 \end{align}
 To verify our claim, notice that, by the very definition of symmetrization with respect to a convex body and by \eqref{2024-61},
 $$(\Tu)^\smallL=u^\smallL =u.$$
 It remains to show that
\begin{equation}\label{Phi_KL}
\left(\Phi_{\bullet\smallK\bullet}\right)_{\bullet\smallL \bullet}=\Phi.
\end{equation}
Thanks to the involution property \eqref{2024-38} of Young conjugation,
equation \eqref{Phi_KL} is equivalent to
\begin{equation}\label{Phi_bu_KL}
\left(\left(\Phi_\bullet\right)_\smallK\right)_\smallL=\Phi_\bullet.
\end{equation}
Owing to Lemma \ref{convesymm}, the assumption \eqref{2024-60} on $\Phi$ implies that $\Phi(-\xi)=A(H^L(\xi))$,
for some  Young function $A$, where $H^L$ is support function of $L$.
Hence, from
\cite[Th. 15.3]{Rockafellar} we deduce  that
\[
\Phi_\bullet(-\xi)=A_\bullet(H_0^L(\xi)),
\]
where $H_0^L$ is the gauge function of $L$. This formula  ensures that the sub-level sets of $\Phi_\bullet$ are homothetic to $-L$, a property which implies \eqref{Phi_bu_KL}. The identity \eqref{Phi_KL}, and hence equation \eqref{2024-63}, are thus established.  Combining \eqref{ineq_klimov_r}, \eqref{ineq_klimov_reverse}, and \eqref{2024-63} yields:
$$\int_{\R^n}\Phi_{\bullet \smallK \bullet}(\nabla \Tu)\, dx = \int_{\R^n} {\Phi}(\nabla u)\, dx.$$
The proof  is complete.
\end{proof}

\begin{proof}[Proof of Proposition~\ref{t:klimov-suff2}]
 The assumption \eqref{growth_infinity} implies that $\Phi_\bullet(\xi)$ is a finite-valued convex function in $\rn$. Hence, any function $u$ having the form \eqref{2024-70prima} is Lipschitz continuous and the support of $\nabla u$ is bounded. As a consequence,  $u\in V^{1,\Phi}_{\rm d}(\R^n)$.
\\ To
 prove that equality holds in \eqref{ineq_klimov}, it suffices to show that the assumptions of Theorem~\ref{t:klimov-suff} are satisfied. Assume $t\in(\essinf u, \esssup u)=(t_1,t_2)$, and define
$$s_t=(t_3-t)/a.$$
 The condition \eqref{2024-30} holds because $s_t>0$.  Since
\[
\{x: u(x)\geq t\}=\bigg\{x: \Phi_\bullet\left(\frac{x_0-x}{a}\right)\leq s_t\bigg\}=\{x_0-ay: \Phi_\bullet(y)\leq s_t\}=-a\{y:\Phi_\bullet(y)\leq s_t\}+x_0,
\]
equation \eqref{t:klimov-suff-a} holds with $a_t=a$ and $x_t=x_0$.
\\ As far as  \eqref{t:klimov-suff-c} is concerned, observe that, for a.e. $t\in(t_1,t_2)$, $\nabla u(x)$ exists for $\cH^{n-1}$-a.e. $x\in \{u=t\}$. Given any $t$ enjoying this property, and any $x\in\{u=t\}$ such that $\nabla u(x)$ exists, one has that
\[
\Phi_\bullet\left(\frac{x_0-x}{a}\right)=s_t
\quad\text{ and } \quad
\nabla u(x)=\nabla \Phi_\bullet\left(\frac{x_0-x}{a}\right).\]
Thus, on setting $\xi=(x_0-x)/a$, we have that $\xi\in \{\Phi_\bullet=s_t\}$ and $\nabla u(x)= \nabla\Phi_\bullet(\xi)$. This proves  \eqref{t:klimov-suff-c}.
\\
A similar argument shows that \eqref{t:klimov-suff-d} holds as well. Indeed,
  $$
 \Tu(x)=\cutoff{t_1}{t_2}\left(t_3-a\Phi_{\bullet \smallK}\left(-\frac{x}{a}\right)\right) \quad \text{for $x\in \rn$.}
$$
Thus,  for a.e. $t\in(t_1,t_2)$, $\nabla \Tu(x)$ exists for $\cH^{n-1}$-a.e. $x\in\{\Tu=t\}$, and, if we define $\xi=-x/a$, then
$\xi\in\{\Phi_{\bullet\smallK}=s_t\}$ and $\nabla \Tu(x)=\nabla\Phi_{\bullet\smallK}(\xi)$.
    \\
Finally, observe that $\nabla u(x)\neq0$ if $u(x)\in(t_1,t_2)$, and $\nabla \Tu(x)\neq0$ if $\Tu(x)\in(t_1,t_2)$. Hence,
$$|\{\nabla u=0, u>t\}|=|\{\nabla \Tu=0, \Tu>t\}|=|\{u=t_2\}|  \quad \text{for $t\in(t_1,t_2)$.}$$
Thereby,
$$\frac{d}{dt}|\{\nabla u=0, u>t\}| = \frac{d}{dt} |\{\nabla \Tu=0, \Tu>t\}|=0 \quad \text{for $t\in(t_1,t_2)$.}$$
Equation \eqref{t:klimov-suff-b} hence follows via \eqref{may1} and \eqref{may2}.
\end{proof}

 \par\noindent {\bf Data availability statement.} Data sharing not applicable to this article as no datasets were generated or analysed during the current study.

\section*{Compliance with Ethical Standards}\label{conflicts}

\smallskip
\par\noindent
{\bf Funding}. This research was partly funded by:
\\ (i) GNAMPA   of the Italian INdAM - National Institute of High Mathematics (grant number not available);
\\ (ii) Research Project   of the Italian Ministry of Education, University and
Research (MIUR) Prin 2017 ``Direct and inverse problems for partial differential equations: theoretical aspects and applications'',
grant number 201758MTR2;
\\ (iii) Research Project   of the Italian Ministry of Education, University and
Research (MIUR) Prin 2022 ``Partial differential equations and related geometric-functional inequalities'',
grant number 20229M52AS, co-funded by PNRR;

\bigskip
\par\noindent
{\bf Conflict of Interest}. The authors declare that they have no conflict of interest.


\begin{thebibliography}
{BGGK}


\bibitem[Al]{alberico}  A.~Alberico, Boundedness of solutions to anisotropic variational problems,  \emph{Comm. Partial Differential Equations} {\bf 36} (2011), 470--486.


\bibitem[ACCZ]{ACCZ}  A.~Alberico, I.~Chlebicka, A.~Cianchi \& A.~Zatorska-Goldstein, Fully anisotropic elliptic problems with minimally integrable data, {\em Calc. Var. Partial Differential Equations } {\bf 58} (2019),  no. 6, Paper No. 186, 50 pp. 35J62.




\bibitem[ADF]{ADF}  A.~Alberico, F.~Feo \& G.~Di Blasio, An eigenvalue problem for the anisotropic $\Phi$-Laplacian, {\em J. Differential Equations } {\bf 269} (2020),  4853--4883.



\bibitem[AFLT]{AFLT}  A.~Alvino, V.~Ferone, P.-L.~Lions \& G.~Trombetti,  Convex symmetrization and applications, {\em Ann. Inst. H. Poincar\'e C Anal. Non Lin\'eaire} {\bf 14} (1997),  275--293.

\bibitem[Bae]{Bae}  A.~Baernstein II,  {\em Symmetrization in analysis},  With David Drasin and Richard S. Laugesen. With a foreword by Walter Hayman.  Cambridge University Press, Cambridge, 2019.

\bibitem[BCF]{BarchiesiCagnettiFusco}  M.~Barchiesi, F.~Cagnetti, \& N.~Fusco,  Stability of the Steiner symmetrization of convex sets, {\em J. Eur. Math. Soc. (JEMS)} {\bf 15} (2013),  1245--1278.

\bibitem[Ba]{barletta}  G.~Barletta,  On a class of fully anisotropic elliptic equations, {\em Nonlinear Anal.} {\bf 197} (2020),  111838, 23 pp.





\bibitem[BGGK]{BiGaGrKi2022}G.~Bianchi, R.~J.~Gardner, P.~Gronchi \& M.~Kiderlen, The P\'olya-Szeg\H{o} inequality for smoothing rearrangements, {\em J. Funct. Anal.} {\bf 287} (2024) 110422, 56 pp.

\bibitem[BrZi]{BrothersZiemer}  J.~E.~Brothers and W.~P.~Ziemer, Minimal rearrangements of Sobolev functions, {\em J. Reine Angew. Math.  (Crelle J.)} {\bf 384} (1988),  153--179.


\bibitem[ChNa]{ChNa} I.~Chlebicka \& P.~Nayar, Essentially fully anisotropic Orlicz functions and uniqueness to measure data problem, {\em Math. Methods Appl. Sci.} {\bf 45} (2022), 8503--8527.

\bibitem[Ci1]{cianchi_pacific}
A.~Cianchi, A fully anisotropic Sobolev inequaliy,
\emph {Pacific J. Math.} {\bf 196} (2000), 283--295.


\bibitem[Ci2]{CiAIHP}
A.~Cianchi,
Local boundedness of minimizers of anisotropic functionals, \emph{Ann. Inst. H. Poincar\'e Anal. Non Lin\'eaire} {\bf 17} (2000), 147--168.



\bibitem[Ci3]{cianchi_ibero}
A.~Cianchi,
Optimal Orlicz-Sobolev embeddings, \emph{Rev. Mat. Iberoamericana} {\bf 20} (2004), 427--474.

\bibitem[Ci4]{CiCPDE}
A.~Cianchi,
Symmetrization in anisotropic elliptic problems,
\emph {Comm. Partial Differential Equations} {\bf 32} (2007), 693--717.


\bibitem[CiFu1]{CianchiFuscoARMA} A.~Cianchi \& N.~Fusco, Functions of bounded variation and rearrangements, {\em Arch. Ration. Mech. Anal.} {\bf 165} (2002),  1–40.

\bibitem[CiFu2]{CianchiFuscoAdvMath} A.~Cianchi \& N.~Fusco, Steiner symmetric extremals in P\'olya-Szeg\"o type inequalities, {\em Adv. Math.} {\bf 203} (2006),  673--728.

\bibitem[CiFu3]{CiFuAppAn} A.~Cianchi \& N.~Fusco, Minimal rearrangements, strict convexity and critical points, {\em Appl. Anal.} {\bf 85} (2006),   67--85.


\bibitem[CiSa]{CianchiSalani} A.~Cianchi \& P.Salani, Overdetermined anisotropic elliptic problems, {\em Math. Ann.} {\bf 345} (2009), 859--881.




\bibitem[EsTr]{EsTr} L.~Esposito \& C.~Trombetti, Convex symmetrization and P\'olya-Szeg\"o inequality, {\em Nonlinear Anal.} {\bf 56} (2004),   43--62.

\bibitem[FeVo]{FeVo} A.~Ferone \& R.~Volpicelli,  Convex rearrangement: equality cases in the P\'olya-Szeg\"o inequality, {\em Calc. Var. Partial Differential Equations} {\bf 21} (2004),   259--272.


\bibitem[FoMu]{FonMul1991} I.~Fonseca \& S.~M\"{u}ller, A uniqueness proof for the Wulff theorem, {\em Proc. Roy. Soc. Edinburgh Sect. A} {\bf 119} (1991), 125--136

\bibitem[Ka]{Ka}
 B.~Kawohl, {\em Rearrangements and convexity of level sets in PDE},   Springer-Verlag, Berlin, 1985.


\bibitem[Ke]{Ke}
 S.~Kesavan, {\em Symmetrization \& applications}, World Scientific Publishing Co. Pte. Ltd., Hackensack, NJ, 2006.

\bibitem[Kl1]{Kli74} V.~S.~Klimov, Isoperimetric inequalities and embedding theorems (Russian), {\em Dokl. Akad. Nauk
SSSR}, {\bf 217} (1974), 272--275; translation {\em Soviet Math. Dokl.}, {\bf 15}.

\bibitem[Kl2]{Kli76} V.~S.~Klimov, Imbedding theorems and geometric inequalities (Russian), {\em Izv. Akad. Nauk SSSR, Ser. Mat.}, {\bf 40} (1976), 645–671; translation {\em Math. USSR Izvestiya}, {\bf 10}, 615--638.


\bibitem[Kl3]{Kli99} V.~S.~Klimov, On the symmetrization of anisotropic integral functionals (Russian), {\em Izv. Vyssh. Uchebn. Zaved. Mat.}, {\bf 8} (1999), 26--32; translation {\em Russian Math. (Iz. VUZ)}, {\bf 43} (1999), 23--29.



\bibitem[LYZ]{Lyz}
E.~Lutwak, D.~Yang \& G.~Zhang, Sharp affine Lp Sobolev inequalities, {\em  J. Differential Geom.} {\bf 62} (2002), 17–38.

\bibitem[Ro]{Rockafellar}
R.~T.~Rockafellar, {\em Convex Analysis}, Princeton University Press, Princeton, N.J., 1970.

\bibitem[Sc]{Sch93}
R.~Schneider, {\em Convex Bodies: The Brunn-Minkowski Theory}, second edition, Cambridge University Press, Cambridge, 2014.

\bibitem[Ta]{Talenti}
G.~Talenti, Best constant in Sobolev inequality, {\em Ann. Mat. Pura Appl.} {\bf 110} (1976), 353--372.

\bibitem[VSch]{VS2006}
J.~Van Schaftingen, Anisotropic symmetrization, {\em Ann. Inst. H. Poincar\'{e} Anal. Non Lin\'{e}aire} {\bf 23} (2006), 539--565.

\bibitem[Wu]{Wu} G.~Wulff, Zur Frage der Geschwindigkeit
des Wachstums und der Aufl\"osung der Kristallfl\"aschen, \emph{Z.
Krist.} {\bf 34} (1901), 449--530.

\bibitem[Zh]{Zhang} G.~Zhang, The affine Sobolev inequality, \emph{J. Differential Geom.} {\bf 53} (1999), 183--202.



\end{thebibliography}
\end{document}